\documentclass[11pt]{amsart}
\usepackage{amscd,amssymb}
\usepackage{amsmath}
\usepackage{graphicx}
\usepackage{subfigure,fancybox}
\usepackage{epic,eepic}
\usepackage{amstext}
\usepackage{amsthm}
\usepackage{amsmath}
\usepackage[utf8]{inputenc}
\usepackage[english]{babel}
\usepackage[square,numbers]{natbib}
\usepackage{booktabs}
\usepackage[hidelinks]{hyperref}
\usepackage[capitalise,nameinlink]{cleveref}
\crefname{equation}{}{}
\crefformat{equation}{\textup{#2(#1)#3}}
\usepackage{comment}
\usepackage{mathtools}
\usepackage{tikz,pgfplots}
\usepackage{tabularx}
\newcolumntype{Y}{>{\centering\arraybackslash}X}
\newcommand{\with}{\,:\,}

\def\N{\mathbb{N}}

\def\R{\mathbb{R}}

\def\calT{\mathcal{T}}

\def\VH{\mathcal{V}_H}
\def\UH{\mathcal{U}_H}
\def\tUH{\tilde{U}_H}

\def\V{\mathcal{V}}
\newcommand{\W}{\mathcal{W}}
\def\tUHl{\tilde{U}_H^\ell}

\def\frakB{\mathfrak{B}}
\def\frakD{\mathfrak{D}}

\def\calR{\mathcal{R}}
\def\calB{\mathcal{B}}
\def\calC{\mathcal{C}}
\def\calCT{\mathcal{C}_T}
\def\calCTl{\mathcal{C}_T_\ell}
\def\calCTl{\mathcal{C}_T^\ell}
\def\calCl{\mathcal{C}^\ell}
\def\IH{\mathcal{I}_H}
\def\EH{\mathcal{E}_H}
\def\PiH{\Pi_H}
\def\PH{\mathcal{P}_H}

\newcommand{\supp}{\operatorname{supp}}

\def\Nb{\mathsf{N}} 

\newcommand{\Ca}{C_\mathrm{a}} 
\newcommand{\Cd}{C_\mathrm{d}} 
\newcommand{\Cbo}{C_\mathrm{bo}} 
\newcommand{\Cb}{C_\mathrm{b}} 
\newcommand{\Ce}{C_\mathrm{e}} 
\newcommand{\Cel}{C_\mathrm{el}} 
\newcommand{\Cl}{C_\mathrm{l}} 
\newcommand{\Ct}{C_\mathrm{t}} 
\newcommand{\Col}{C_\mathrm{ol}} 
\newcommand{\Cp}{C_{\mathcal P_H}} 
\newcommand{\CI}{C_{\mathcal I_H}} 
\newcommand{\CF}{C_\mathrm{PF}} 

\def\lamjT{\Lambda_{j,T}}

\def\bjT{b_{j,T}}
\def\vH{v_H}

\newcommand{\lamvH}{\lambda_{\vH}}
\newcommand{\muH}{\mu_H}

\newcommand{\qH}{q_H}

\newcommand{\vC}{\tilde{v}_H}

\newcommand{\uC}{\tilde{u}_H}

\newcommand{\uCl}{\tilde{u}_H^\ell}
\newcommand{\vCl}{\tilde{v}_H^\ell}

\newcommand\dx{\,\text{d}x}

\newcommand{\tnorm}[2]{{\left\lVert #1 \right\rVert}_{L^2(#2)}}

\newcommand{\tnormf}[2]{{\| #1 \|}_{L^2(#2)}}

\newcommand{\vnorm}[2]{{\left\lVert #1 \right\rVert}_{\V(#2)}}

\newcommand{\tspf}[3]{{( #1\,,\,#2 )}_{L^2(#3)}}

\newtheorem{theorem}{Theorem}[section]

\newtheorem{lemma}[theorem]{Lemma}

\theoremstyle{definition}

\theoremstyle{remark}
\newtheorem*{remark}{Remark}
\numberwithin{theorem}{section}
\numberwithin{equation}{section}
\numberwithin{table}{section}
\numberwithin{figure}{section}

\allowdisplaybreaks

\begin{document}
	
\title[An improved high-order multiscale method]{An improved high-order method for elliptic multiscale problems}
\author[Z.~Dong, M.~Hauck, R.~Maier]{Zhaonan Dong$^\ast$, Moritz Hauck$^\dagger$, Roland Maier$^\ddagger$}
\address{${}^{\ast}$ Inria, 2 rue Simone Iff, 75589 Paris, France \&  CERMICS, Ecole des Ponts, 77455 Marne-la-Vall\'{e}e 2, France}
\email{zhaonan.dong@inria.fr}
\address{${}^{\dagger}$ Institute of Mathematics, University of Augsburg, Universit\"atsstr.~12a, 86159 Augsburg, Germany}
\email{moritz.hauck@uni-a.de}
\address{${}^{\ddagger}$ Institute of Mathematics, Friedrich Schiller University Jena, Ernst-Abbe-Platz 2, 07743 Jena, Germany}
\email{roland.maier@uni-jena.de}

\thanks{The work of Moritz Hauck is part of a project that has received funding from the European Research Council (ERC) under the European Union's Horizon 2020 research and innovation programme (Grant agreement No.~865751 --  RandomMultiScales).}

\begin{abstract}
In this work, we propose a high-order multiscale method for an elliptic model problem with rough and possibly highly oscillatory coefficients. Convergence rates of higher order are obtained using the regularity of the right-hand side only. Hence, no restrictive assumptions on the coefficient, the domain, or the exact solution are required. 
In the spirit of the Localized Orthogonal Decomposition, the method constructs coarse problem-adapted ansatz spaces by solving auxiliary problems on local subdomains. More precisely, our approach is based on the strategy presented by Maier [SIAM J.~Numer.~Anal.~59(2),~2021]. The unique selling point of the proposed method is an improved localization strategy curing the effect of deteriorating errors with respect to the mesh size when the local subdomains are not large enough.  We present a rigorous a priori error analysis and demonstrate the performance of the method in a series of numerical experiments. 
\end{abstract}
\maketitle
\indent
{\tiny {\bf AMS subject classifications.}  
	{\bf 65N12}, 
	{\bf 65N30} 
} 
%
\section{Introduction}

In this work, we consider the numerical solution of elliptic partial differential equations (PDEs) with possibly rough and highly varying coefficients. The numerical treatment of such problems with classical finite element approaches suffers from suboptimal approximation rates and pre-asymptotic effects if the discretization parameters do not resolve the oscillation scales. As globally resolving these scales can be computationally unfeasible, we aim to obtain reasonable errors already on coarse scales. Computational multiscale methods tackle this problem by constructing coarse problem-adapted ansatz spaces that contain local information on the underlying coefficient. Under minimal regularity assumptions on the coefficient, such methods achieve optimal orders of approximation while introducing only a moderate computational overhead (e.g., an increased support of the basis functions or an increased number of basis functions per mesh entity). Prominent first-order methods include, e.g., Generalized (Multiscale) Finite Element Methods~\cite{BabO83,BabCO94,BabL11,EfeGH13}, Adaptive Local Bases~\cite{GraGS12,Wey16}, the Localized Orthogonal Decomposition (LOD) method~\cite{MalP14,HenP13}, or gamblets~\cite{Owh17}. For an overview of such methods, see also the recent textbooks~\cite{OwhS19,MalP20} and the review article~\cite{AltHP21}.

Under additional smoothness assumptions, higher-order multiscale methods have been proposed based on, e.g., the Heterogeneous Multiscale Method~\cite{LiPT12,AbdB12} or the Multiscale Finite Element Method~\cite{AllB05,HesZZ14}. Other strategies include Multiscale Hybrid-Mixed Methods~\cite{HarPV13,AraHPV13} as well as the Multiscale Hybrid High-Order Method~\cite{CicEL19}. However, for general $L^\infty$-coefficients, these methods only provide first-order convergence results in $H^1$. This issue is overcome in~\cite{Mai21} (see also~\cite{Mai20}), where a method is proposed that achieves high-order rates for arbitrarily rough coefficients. The idea is to exploit appropriate orthogonality properties in the spirit of the LOD which allows one to extract higher-order rates from the right-hand side only. Basis functions of the coarse problem-adapted ansatz space can be constructed as solutions to certain constrained energy minimization problems. The constraints enforce that the $L^2$-projection of a particular basis function into the space of piecewise polynomials coincides with a tensor-product Legendre polynomial of a given degree that is only supported on a single (quadrilateral) element. Since these basis functions are globally supported, a localization of the basis is required in order to obtain a practically feasible approach. In~\cite{Mai21}, it is shown that the problem-adapted basis functions decay exponentially fast in modulus relative to the coarse mesh, which allows one to approximate these functions on localized subdomains. However, the approach theoretically and practically suffers from deteriorating errors as the mesh size is decreased when the local subdomains are not suitably increased.

In this work, we cure the above-mentioned error deterioration of the method in~\cite{Mai21} by proposing an improved localization technique. It slightly enlarges the support of the lowest-order basis functions, while keeping the other basis functions unchanged. Notably, this cures the issue of deteriorating errors for {an arbitrary but fixed} polynomial degree. In practice, this modification avoids an overestimation of the required size of the local subdomains and thereby reduces the computational costs and increases the reliability of the method.

This work is structured as follows. First, we present the elliptic model problem (\Cref{s:model}) and preliminary notation (\Cref{s:prelim}). A prototypical (non-localized) multiscale method is then presented in  \Cref{s:idealmethod}. In contrast to the definition in~\cite{Mai21}, which is based on non-conforming polynomial spaces, our construction utilizes conforming bubble spaces. These spaces are designed such that the resulting prototypical multiscale method concides with the one in~\cite{Mai21}. In \Cref{s:decloc}, we propose an improved localization technique, which allows one to derive a practical version of the method (\Cref{s:practical}) which is able to remedy the issues of deteriorating errors for subdomains of fixed size.
In \Cref{s:numerics}, we finally show the performance of the novel strategy in a set of numerical experiments and, in particular, compare it to the original version.

\section{Model problem}\label{s:model}
On a polygonal Lipschitz domain $D \subset \R^d$, $d \in \N$, we consider the prototypical second-order elliptic problem in weak form that seeks $u \in H^1_0(D)$ such that
\begin{equation}\label{eq:prob}
	a(u, v)\coloneqq  \int_D A \nabla u \cdot \nabla v \dx = (f, v)_{L^2(D)} \qquad\text{for all }v \in H^1_0(D)
\end{equation}
with a right-hand side $f \in L^2(D)$ and a scalar-valued coefficient $A \in L^\infty(D)$ that fulfills $0 < \alpha \leq A(x) \leq \beta < \infty$ for almost all $x \in D$ and some positive constants $\alpha,\beta$. 
Apart from these bounds, we do not pose any further assumptions on $A$. 
However, we implicitly assume that the coefficient oscillates on multiple scales. 
It is well-known that in such a setting, classical finite element methods can perform arbitrarily badly~\cite{BabO00}, in particular if the oscillations in the coefficient are not resolved.
In the context of the model problem~\cref{eq:prob}, this is also illustrated in~\cite[Ch.~2]{MalP20}.
We emphasize that the choice of scalar-valued coefficients is mainly to simplify the presentation. All results also hold for matrix-valued symmetric coefficient functions. In that case, $\alpha$ and $\beta$ denote the uniform lower and upper spectral bounds, respectively.  

By the Lax--Milgram theorem, the solution $u$ exists and is uniquely defined. Further, the stability estimate
\begin{equation}\label{eq:stabu}
	\|\nabla u \|_{L^2(D)} \leq \CF\alpha^{-1} \|f\|_{L^2(D)}
\end{equation}
holds, where the Poincar\'e--Friedrichs inequality with constant $\CF>0$ is used to bound the $H^{-1}(D)$-norm of the right-hand side by its $L^2(D)$-norm.

\section{Preliminaries}\label{s:prelim}

Before we state a prototypical multiscale method in the following section, we first need to introduce coarse-scale spaces and interpolation operators that are required for the construction of the method.

\subsection{Non-conforming spaces}
Let $\calT_H$ be a quasi-uniform and non-degenerate quadrilateral mesh with mesh size $H > 0$. We define for $\ell \in \N$ the \emph{element patch of order~$\ell$} around $S \subset D$ by
\begin{equation}
	\label{eq:patch}
	\Nb^\ell(S) \coloneqq  \Nb^1(\Nb^{\ell-1}(S)),\quad \ell \geq 2,\qquad\Nb^1(S) \coloneqq  \bigcup \bigl\{T \in \calT_H\,\colon\, \overline{S} \,\cap\, \overline{T}\neq \emptyset\bigl\}
\end{equation}
and set $\Nb(S) \coloneqq \Nb^1(S)$.
For a fixed (but arbitrary) polynomial degree $p$, we denote with
\begin{equation*}
	\VH \coloneqq  \{v \in L^2(D)\;\colon\; v\vert_T, T \in \calT_H, \text{ is a polynomial of coordinate degree }\leq p\}
\end{equation*}
the non-conforming space (with respect to $H^1_0(D)$) consisting of element-wise defined polynomials.
The restriction of $\VH$ to a subdomain $S \subset D$ is given by
\begin{equation*}
	\VH(S) \coloneqq  \{v \in \VH\;\colon\; \supp(v) \subset S\}.
\end{equation*}
As orthonormal basis $\frakD_T \coloneqq  \{\lamjT\}_{j=1}^{N}$ of the space $\VH(T)$, $T \in \calT_H$, we henceforth use the shifted tensor-product Legendre polynomials. Thus, $\frakD \coloneqq  \bigcup_{T \in \calT_H} \frakD_T$ defines a local basis of $\VH$. Denoting with $M \coloneqq  |\calT_H|$ the number of elements in $\calT_H$, this means that $\dim \VH = MN$. Note that it holds $N = (p+1)^d$.

Let $\PiH\colon L^2(D) \to \VH$ be the $L^2$-projection defined for any $v \in L^2(D)$ by the element-wise equation
\begin{equation*}\label{eq:L2proj}
	(\PiH v, \vH)_{L^2(T)} = (v,\vH)_{L^2(T)} \qquad\text{for all }\vH \in \VH(T),\;T \in \calT_H.
\end{equation*}
From this definition, we directly obtain the following stability estimate,
\begin{equation*}\label{eq:L2stab}
	\|\PiH v\|_{L^2(T)} \leq \|v\|_{L^2(T)},\qquad v \in L^2(T),\;T \in \calT_H.
\end{equation*}
Additionally, there exists a constant $\Ca > 0$ (independent of the mesh size) such that
\begin{equation}\label{eq:approx}
	\|(1-\PiH) v\|_{L^2(T)} \leq \Ca H \|\nabla v\|_{L^2(T)},\qquad v \in H^1(T),\;T \in \calT_H,
\end{equation}
see, e.g., \cite{HouSS02}.
We emphasize that the constant $\Ca$ depends on the polynomial degree~$p$. {More precisely, $\Ca \sim (p+1)^{-1}$.} 
Note that the above definitions are based on quadrilateral meshes. However, this is not necessary and all results can be derived analogously for simplicial meshes after suitably adapting the definitions. For instance, $\VH$ would be defined as the space of piecewise polynomials of total degree $p$ (instead of coordinate degree) and the number of degrees of freedoms changes accordingly.

For late use, we also define the \emph{broken Sobolev space} $H^k(\calT_H)$, $k \in \N$, by
\begin{equation*}
	H^k(\calT_H) \coloneqq \big\{ v \in L^2(D)\;:\; \forall T \in \calT_H\;:\; v\vert_T \in H^k(T) \big\}
\end{equation*}
with the seminorm
\begin{equation*}
	|\cdot|^2_{H^k(\calT_H)} \coloneqq \sum_{T \in \calT_H} |\cdot|^2_{H^k(T)}.
\end{equation*}

\subsection{Conforming spaces}
While the locality of the non-conforming spaces is very beneficial for designing multiscale methods,
it makes their theoretical investigation more cumbersome.
Therefore, we replace the space $\VH$ by a conforming space $\UH \subset H^1_0(D)$ with the properties that $\PiH \UH = \VH$ and $\dim \UH = \dim \VH = MN$.
The space $\UH$ is constructed based on a set of local basis \emph{bubble functions} $\frakB_T \coloneqq  \{\bjT\}_{j=1}^N$ for each element $T \in \calT_H$. These functions fulfill, for $1 \leq j \leq N$,
\begin{equation}\label{eq:bubble}
	\bjT \in H^1_0(T), \qquad \PiH \bjT = \lamjT.
\end{equation}
Note that, throughout this paper, we do not distinguish between $H^1_0(S)$-functions on a subdomain
$S \subset D$ and
their $H^1_0(D)$-conforming extension by zero to the full domain $D$.

For the moment, an explicit characterization of the basis
$\frakB \coloneqq  \bigcup_{T \in \calT_H} \frakB_T$ of the space $\UH$ is not required. However, it is important that such functions actually exist, which is stated in the following lemma.

\begin{lemma}[Local bubble functions]\label{lem:locbubble}
	There exists a constant $\Cb>0$ which solely depends on the regularity of $\calT_H$ and the polynomial degree $p$, such that for all $T \in \calT_H$ and $j = 1,\dots,N$, there exists a function $\bjT \in H^1_0(T)$ with $\PiH \bjT  = \lamjT$ and
	\begin{equation*}
		\Big\|\nabla {\sum_{j=1}^N c_j}\, \bjT\Big\|_{L^2(T)} \leq \Cb H^{-1}\, \Big\|{\sum_{j=1}^N c_j}\,\lamjT\Big\|_{L^2(T)}
	\end{equation*}
	{for any coefficients $c_j \in \R$. Regarding the scaling of $\Cb$ with respect to $p$, we have $\Cb \sim  (p+1)^2$.}
\end{lemma}

\begin{proof}
	The proof is a special case of~\cite[Cor.~3.6]{Mai21}.
\end{proof}

Based on the bases $\frakD$ of $\VH$ and $\frakB$ of $\UH$, we can construct a continuous mapping $\calB_H\colon \VH \to \UH \subset H^1_0(D)$ that is uniquely defined by the property that
\begin{equation*}
	\calB_H \lamjT = \bjT\qquad\text{ for all }T \in \calT_H,\,1\leq j\leq N.
\end{equation*}
This \emph{bubble operator} $\calB_H$ can be extended to functions in $L^2(D)$ with the definition
\begin{equation*}
	\calB_H v \coloneqq \calB_H \PiH v\qquad\text{for all } v \in L^2(D)
\end{equation*}
and it can be shown that it has the following representation,
\begin{equation}
	\calB_H v = \sum_{T \in \calT_H}\sum_{j = 1}^{N} \Big(\int_T v \lamjT \dx\Big) \bjT,\qquad v \in L^2(D).\label{eq:repB}
\end{equation}
It is straightforward that $\ker \PiH = \ker \calB_H$. {Using  \Cref{lem:locbubble}, we directly have}
\begin{equation}\label{eq:inv}
	\|\nabla \calB_H v\|_{L^2(T)} \leq \Cbo H^{-1} \|v\|_{L^2(T)},\qquad v \in L^2(T),\;T \in \calT_H.
\end{equation}
The operator $\calB_H$ will be important for the definition of a prototypical multiscale method in the following section and to link it to the method stated in~\cite{Mai21}.

For later use, we introduce a quasi-interpolation operator $\PH\colon L^2(D) \to H^1_0(D)$ which has the same kernel as the $L^2$-projection $\PiH$. The construction of $\PH$ is based on another (low-order) quasi-interpolation operator defined as $\IH\coloneqq \EH\circ \PiH^0$, where $\PiH^0$ denotes the $L^2$-projection onto piecewise constants and $\EH$ is the averaging operator mapping piecewise constants to the space of continuous piecewise linear functions. For inner nodes~$z$, we define
\begin{equation*}
	(\EH v)(z) \coloneqq  \frac{1}{\mathrm{card}\{T \in \calT_H\,\colon\, z \in \overline{T}\}}\sum_{T\in\calT_H\,:\, z \in \overline{T}} v\vert_T(z)
\end{equation*}
and we set $(\EH v)(z) \coloneqq 0$ for nodes at the boundary. The operator $\EH$ is well-known from the theory of domain decomposition methods, see \cite{Osw93,Bre94}. 
Using the bubbles $b_{j,T}$, we then locally manipulate the operator $\IH$ such that its kernel coincides with the kernel of $\PiH$. More precisely, for any $v \in L^2(D)$ we define
\begin{equation}\label{eq:defPH}
	\PH v \coloneqq \IH v + \calB_H(v -\IH v).
\end{equation}
The operator $\PH$ has been introduced in~\cite[eq.~(3.12)]{AltHP21} employing very general functionals, so-called \emph{quantities of interest}.
{
	\begin{remark}[Effect of $\PH$]
		In the definition of $\PH$, the additional term involving $\IH$ only has an effect on the zero-order bubble function (cf.~\Cref{lem:applicationPH} below for a proof). That is, the image of $\PH$ is given by the span of the basis functions $\frakB$ of $\UH$, except for the zero-order bubble function, which is replaced by a modified function with a slightly increased support. The motivation for this modification is that $\PH$ is better suited than the operator $\calB_H$ to approximate (piecewise) constant functions. This property is crucial for an $H$-independent stability estimate as proved in \Cref{lem:propPH} and will be important to obtain an improved version of the method in~\cite{Mai21}.
	\end{remark}
}
The precise properties of $\PH$ are investigated in the following lemma.

\begin{lemma}[Properties of $\PH$]\label{lem:propPH}
	The operator $\PH$ is a projection, i.e., it satisfies $\PH^2 = \PH$. Moreover, its kernel coincides with the kernel of $\PiH$, i.e.,
	\begin{equation*}
		\ker \PiH = \ker\PH.
	\end{equation*}
	Further, there exists a constant $\Cp>0$, which solely depends on the regularity of $\calT_H$ and the polynomial degree $p$, {more precisely, $\Cp \sim (p+1)^{2}$,} such that for all $v \in H^1_0(D)\vert_{\Nb(T)},\;T \in \calT_H$
	\begin{align}
		\tnormf{\nabla \PH v}{T}  +  H^{-1}\tnormf{(1-\PH)v}{T} \leq\Cp \tnormf{ \nabla v}{\mathsf N(T)}\label{eq:PHprop}.
	\end{align}
\end{lemma}

\begin{proof}
	We first prove that $\PH$ is a projection.  A technical calculation shows that $\PH$ admits the representation
	\begin{equation}
		\label{eq:repPH}
		\PH v = \sum_{T\in \calT_H}\sum_{j=1}^N \Big(\int_T v \lamjT\dx\Big) \PH b_{j,T}.
	\end{equation}
	Using this representation, the projection property follows from the identity
	\begin{equation*}
		\int_T \PH b_{i,K} \Lambda_{j,T}\dx = \delta_{ij}\delta_{KT},
	\end{equation*}
	where $\delta$ denotes the Kronecker delta function.
	
	The equality of the kernels follows directly form~\cref{eq:repPH} using the fact that the values $\int_Tv \lamjT\dx$ are also the coefficients of $\PiH v$ when expanded with respect to the basis function $\lamjT$.
	
	For the proof of~\cref{eq:PHprop}, we use the stability and approximation properties of $\IH$, see, e.g.,~\cite{ErnG17}. There exists a constant $\CI>0$ such that for all $v \in H^1_0(D)\vert_{\Nb(T)}$ and $T \in \calT_H$,
	\begin{align}
		\tnormf{\nabla \IH v}{T}  +  H^{-1}\tnormf{(1-\IH)v}{T} \leq\CI \tnormf{ \nabla v}{\mathsf N(T)}\label{eq:IHprop}.
	\end{align}
	Using the stability estimate of the bubble operator~\cref{eq:inv} and the properties of~$\IH$ stated in~\cref{eq:IHprop}, we obtain
	\begin{align*}
		\|\nabla \PH v\|_{L^2(T)}&\leq  \|\nabla \IH v\|_{L^2(T)} + \Cbo H^{-1}\|v -\IH v\|_{L^2(T)}\\&\leq \CI(1+\Cbo) \|\nabla v\|_{L^2(\Nb(T))},
	\end{align*}
	which proves the $H^1$-stability of $\PH$ in~\cref{eq:PHprop}. {The scaling of~$\Cp$ with respect to~$p$ follows directly by the corresponding estimate for $\Cbo$.} The second ingredient, the $L^2$-ap\-prox\-i\-ma\-tion property in~\cref{eq:PHprop}, follows with similar arguments and is omitted for the sake of brevity.
\end{proof}

\section{A prototypical multiscale method}\label{s:idealmethod}

The conformity of the space $\UH$ is essential for the construction of a prototypical multiscale method, which directly relates to the classical LOD methodology as presented in~\cite{MalP14,HenP13,MalP20}. We will adapt ideas from~\cite{HauP21} and prove that the constructed $p$th-order prototypical multiscale method is equivalent to the one proposed in~\cite{Mai21}. The benefit of this alternative construction is that it enables an improved localization compared to~\cite{Mai21} that will be shown in the subsequent section based on the operator~$\PH$. {To provide a better overview of the important spaces and operators, we have provided a summary of the most important definitions in \Cref{table: spaces and operators} below.}

In the spirit of the LOD, we define a so-called fine-scale space $\W$ of functions that are not captured by the $L^2$-projection $\PiH$,
\begin{equation*}
	\W \coloneqq  \ker \PiH\vert_{H^1_0(D)}.
\end{equation*}
Note that the spaces $\UH$ and $\W$ define an $L^2$-orthogonal decomposition of the space $H^1_0(D)$, that is
\begin{equation*}
	H^1_0(D) = \UH \oplus \W,\qquad\quad H^1_0(D) \ni v = \calB_H v + (1-\calB_H) v .
\end{equation*}
However, this decomposition is not suited for the approximation of the solution to~\cref{eq:prob} and we would like to modify the space $\UH$  such that the decomposition is $a$-orthogonal instead and thus includes information on the problem at hand.
Based on the space $\W$, we therefore introduce a so-called \emph{correction operator} $\calC\colon H^1_0(D) \to \W$. For any $v \in H^1_0(D)$, the correction $\calC v \in \W$ is defined as the solution to
\begin{equation}\label{eq:corr}
	a(\calC v,w) = a(v,w) \qquad\text{for all } w \in \W.
\end{equation}
By construction, this defines an $a$-orthogonal decomposition of $H^1_0(D)$ into two subspaces. More precisely, we set $\tUH \coloneqq  (1-\calC)H^1_0(D)$ and observe that
\begin{equation*}
	H^1_0(D) = \tUH \oplus \W,\qquad\quad H^1_0(D) \ni v = (1-\calC) v + \calC v.
\end{equation*}
Note that for the definition of $\tUH$, it suffices to consider the coarse functions $\UH$ in the domain of $(1-\calC)$.
That is, with the above operators, the property $\PiH(1-\calB_H)v = 0$ for all $v \in H^1_0(D)$, and the fact that $\calC$ is the identity on $\W$, we deduce
\begin{equation}\label{eq:deftUH}
	\begin{aligned}
		\tUH &= (1-\calC)H^1_0(D) = (1-\calC)(1-\calB_H)H^1_0(D) + (1-\calC)\calB_H H^1_0(D) \\&= (1-\calC)\calB_H H^1_0(D) = (1-\calC)\calB_H\VH = (1-\calC)\UH.
	\end{aligned}
\end{equation}
Therefore, $(1-\calC)\colon\UH \to \tUH$ is a bijective operator with inverse $\calB_H\vert_{\tUH}$ and it holds $\dim \UH = \dim\tUH$.
Further,
\begin{equation}
	\label{eq:idealbasis}
	\tilde\frakB\coloneqq \bigcup_{b \in \frakB}(1-\calC)b
\end{equation}
defines a basis of $\tUH$.

The newly defined space $\tUH$ is a suitable candidate for approximating the solution $u$ to~\cref{eq:prob} and we define a Galerkin method based on the space $\tUH$ as follows: find $\uC \in \tUH$ such that
\begin{equation}\label{eq:msprob}
	a(\uC,\vC) = (f,\vC)_{L^2(D)}\qquad\text{for all } \vC \in \tUH.
\end{equation}
The well-posedness of~\cref{eq:msprob} follows directly from the Lax-Milgram theorem with a stability estimate as in~\cref{eq:stabu}. The method stated in~\cref{eq:msprob} is typically called the \emph{prototypical method}, because it is defined based on the global correction operator $\calC$ that needs to be suitably localized and discretized in practice. This will be investigated in the following sections.

Before we turn to an optimal error estimate for the solution $\uC$ to~\cref{eq:msprob}, we show that the approximation space $\tUH$ coincides with the multiscale space $\calR\VH$ from
~\cite[Sec.~3.1]{Mai21}. Therein, the operator $\calR$ is defined by seeking for any $\vH \in \VH$ a pair $(\calR\vH,\lamvH) \in H^1_0(D) \times \VH$ solving
\begin{equation}\label{eq:saddlepoint}
	\begin{aligned}
		a(\calR \vH, v)\qquad\, &+& (\lamvH, v)_{L^2(D)} \quad&=\quad 0,\\
		(\calR \vH,\muH)_{L^2(D)} && &=\quad (\vH, \muH)_{L^2(D)}
	\end{aligned}
\end{equation}
for all $v \in H^1_0(D)$ and all $\muH \in \VH$. Note that the unique solvability of~\cref{eq:saddlepoint} is shown in~\cite{Mai21}.
\begin{table}[!t]
	{
		\caption{Important spaces and operators} \label{table: spaces and operators}
		\begin{tabular}{|c|l|}
			\hline
			$\VH$  &   piecewise polynomial space with coordinate degree less than or equal to $p$   \\
			\hline
			$\UH$  &   $H^1_0(D)$-conforming space  with $\PiH \UH = \VH$ and $\dim \UH = \dim \VH$ \\
			\hline
			$\PiH$ &  $L^2$-projection onto $\VH$ \\
			\hline
			$\IH$ &  first-order conforming quasi-interpolation operator  \\
			\hline
			$\calB_H$&  bubble operator $\calB_H\colon \VH \to \UH \subset H^1_0(D)$.  \\
			\hline
			$\PH$ & $\PH \bullet \coloneqq \IH \bullet + \calB_H(\bullet -\IH \bullet)$ with $\ker \PH = \ker\PiH$   \\
			\hline
			$\W$ &  fine-scale space $\W\coloneqq  \ker \PiH\vert_{H^1_0(D)}$\\
			\hline
			$\calC$ &   correction operator $\calC\colon H^1_0(D) \to \W$  \\
			\hline
			$\tUH$ &  $a$-orthogonal complement of $\W$ in $H^1_0(D)$,  $\tUH =  (1-\calC)\UH$\\
			\hline
	\end{tabular}}
\end{table}

The equivalence of spaces is stated in the following lemma.
\begin{lemma}[Equivalence of spaces]\label{lem:equalSpaces}
	It holds $\calR\VH = \tUH$. More precisely, for any $\vH \in \VH$ we have the identity
	\begin{equation*}
		\calR\vH = (1-\calC)\calB_H\vH.
	\end{equation*}
\end{lemma}

\begin{proof}
	Let $\vH \in \VH$. We show that $(1 - \calC)\calB_H\vH$ solves~\cref{eq:saddlepoint} together with an appropriate Lagrange multiplier $\lambda_{\vH}$. The assertion then follows by the uniqueness of~\cref{eq:saddlepoint}.
	
	By~\cref{eq:corr}, we have that $a((1-\calC)\calB_H\vH,w) = 0$ for all $w \in \W$. Let now $\muH \in \VH$ and $v \in H^1_0(D)$. Due to the $L^2$-orthogonality of the spaces $\VH$ and $\W$, we obtain
	\begin{equation}\label{eq:proof_altform_1}
		((1-\calC)\calB_H\vH,\muH)_{L^2(D)} = (\calB_H\vH,\muH)_{L^2(D)} = (\vH,\muH)_{L^2(D)},
	\end{equation}
	which shows that $(1-\calC)\calB_H\vH$ solves the second equation of~\cref{eq:saddlepoint}. We now turn to the first equality. Since $(1-\calB_H)v \in \W$, we have that
	\begin{equation}\label{eq:proof_altform_2}
		a((1-\calC)\calB_H\vH,v) = a((1-\calC)\calB_H\vH,\calB_H v).
	\end{equation}
	To fulfill~\cref{eq:saddlepoint}, we define a suitable Lagrange multiplier. Let $\lambda_{\vH} \in \VH$ be the solution to
	\begin{equation}\label{eq:proof_altform_3}
		(\lambda_{\vH},\qH)_{L^2(D)} = - a((1-\calC)\calB_H\vH,\calB_H\qH)\qquad\text{for all }\qH \in \VH,
	\end{equation}
	which is well-defined due to the coercivity of $(\cdot,\cdot)_{L^2(D)}$ with respect to $\VH \subset L^2(D)$ and the boundedness of $a$ and $\calB_H$.
	The combination of \cref{eq:proof_altform_2} and \cref{eq:proof_altform_3} leads to
	\begin{equation}\label{eq:proof_altform_4}
		a((1-\calC)\calB_H\vH,v) + (\lambda_{\vH},v)_{L^2(D)} = 0,
	\end{equation}
	which concludes the proof.
\end{proof}

\subsection{Error estimates}

The equality of the space $\calR\VH$ and $\tUH$ implies that the error estimates of the prototypical method~\cref{eq:msprob} follow directly from the corresponding ones in~\cite[Sec.~3.2]{Mai21}. In particular, we have the following theorem.
\begin{theorem}[Error of the prototypical method]\label{t:errODp}
	Assume that $f \in H^k(\calT_H)$, $k \in \N_0$, and set $s \coloneqq  \min\{k,p+1\}$. Further, let $u \in H^1_0(D)$ and $\uC \in \tUH$ be the solutions of \cref{eq:prob} and \cref{eq:msprob}, respectively. Then there exists a constant $\Ce>0$, which solely depends on the regularity of $\calT_H$ and the polynomial degree $p$, such that
	\begin{align}
		\|\nabla(u - \uC)\|_{L^2(D)} &\leq \Ce\,H^{s+1}\, |f|_{H^s(\calT_H)},\label{eq:errODp}\\
		\|u - \uC\|_{L^2(D)} &\leq \Ce\,H^{s+2}\, |f|_{H^s(\calT_H)} \label{eq:errODpL2}
	\end{align}
	with the notation $H^0(\calT_H) \coloneqq  L^2(D)$ and $|\cdot|_{H^0(\calT_H)}\coloneqq  \|\cdot\|_{L^2(D)}$.
	{Note that the constant~$\Ce$ generally depends in a positive way on $p$ and $s$ and scales (at least) like~$(p+1)^{-1}$.}
\end{theorem}
\begin{proof}
	The result is a direct consequence of~\cite[Thm.~3.1]{Mai21}.
\end{proof}

\section{Decay and localization}\label{s:decloc}

As mentioned above, the prototypical method defined in~\cref{eq:msprob} is not a practically feasible method since the corresponding basis functions are globally defined. In this section, we show that the Green's function of the correction decays exponentially fast, which justifies the localization of the basis functions. For enhanced stability properties, this paper employs a localization technique based on localized element corrections.

\subsection{Exponential decay}

For any $T \in \calT_H$, the \emph{element correction operator} $\calCT\colon H^1_0(D) \to \W$ is defined, for any $v \in H^1_0(D)$, by the solution $\calCT v \in \W$ to
\begin{align}\label{eq:elemcorrdef}
	a(\calCT v,w) &= a_T(v,w) \qquad \text{for all }w\in \W,
\end{align}
where we use the restricted bilinear form $a_T$ defined by
\begin{equation*}
	a_T(v,w) \coloneqq  \int_TA\nabla v\cdot \nabla w \dx, \qquad v,w \in H^1_0(D).
\end{equation*}
We can then split the correction operator $\calC$ into its element contributions,
\begin{equation}
	\label{eq:sumcorr}
	\calC = \sum_{T\in \calT_H}\calCT.
\end{equation}
The following lemma proves that the moduli of the element corrections decay exponentially fast relative to the coarse mesh.

\begin{lemma}[Exponential decay of element corrections]\label{lemma:expdec}
	There exists a constant $\Cd>0$, which is independent of $H,\ell,T$, such that for all $T \in \calT_H$, $v \in H^1_0(D)$, $\ell \in \N$,
	\begin{equation*}
		\tnormf{\nabla \calCT v}{D\backslash \mathsf N^\ell(T)}\leq \exp(-\Cd\ell )\tnorm{\nabla \calCT v}{D}.
	\end{equation*}
	{Regarding the polynomial degree $p$, we have $\Cd \sim (p+1)^{-1}$.}
\end{lemma}

\begin{proof}
	The proof is presented in \Cref{sec:appendix}.
\end{proof}

{We emphasize that the scaling of the constant $\Cd$ in \Cref{lemma:expdec} is similar to the scaling in~\cite[Thm.~4.1]{Mai21}. However, as pointed out in~\cite[Rem.~4.3]{Mai21}, the negative scaling with respect to $p$ seems to be pessimistic and is possibly not sharp. In fact, a positive scaling with respect to $p$ is observed in our numerical experiments, see~\Cref{fig:dec}.}

\subsection{Localized corrections}
The exponential decay of the element corrections $\calCT$ motivates a localization to element patches around the element $T$. For a given oversampling parameter $\ell\in \N$, we therefore define a localized fine-scale space as
\begin{equation*}
	\W_T^\ell \coloneqq  \W(\mathsf N^\ell(T))\coloneqq \{w\in \W\,:\, \supp(w)\subset \mathsf N^\ell(T)\}\subset \W.
\end{equation*}
For any $T \in \calT_H$, the \emph{localized element correction operator} $\calCTl\colon H^1_0(D) \to \W_T^\ell$ is defined, for any $v \in H^1_0(D)$, by the solution $\calCTl v \in \W_T^\ell$ to
\begin{align}\label{eq:locelemcorr}
	a(\calCTl v,w) = a_T(v,w) \qquad \text{for all }w\in \W_T^\ell.
\end{align}
Without localization, the correction operator $\calC$ equals the sum of the elements contributions $\calCT$; see \cref{eq:sumcorr}. Thus, we also define the localized correction operator $\calCl$ as the sum of the localized element contributions,
\begin{equation*}
	\calCl \coloneqq  \sum_{T\in \calT_H}\calCTl.
\end{equation*}
The next lemma shows exponential approximation properties of this localized correction operator.
\begin{lemma}[Localization error of correctors]\label{lemma:locerrcorr}
	There exists a constant $\Cl>0$ independent of $H,\ell$ such that for all $v \in H^1_0(D)$, $\ell \in \N$
	\begin{equation*}
		\tnormf{\nabla(\calC-\calCl)v}{D}\leq \Cl\ell^{d/2} \exp(-\Cd\ell) \tnormf{\nabla v}{D}
	\end{equation*}
	with the constant $\Cd$ from \Cref{lemma:expdec}. {With respect to the polynomial degree, we have the scaling $\Cl \sim (p+1)$.}
\end{lemma}

\begin{proof}
	The proof is stated in \Cref{sec:appendix}.
\end{proof}

\section{Practical multiscale method}\label{s:practical}

In this section, we present a practical (localized) multiscale method based on the considerations in the previous sections. As mentioned above, the (prototypical) basis $\tilde\frakB$ defined in~\cref{eq:idealbasis} consists of non-local but rapidly decaying basis functions; see \Cref{lemma:expdec}. In order to avoid unfeasible global computations, we replace the prototypical basis by its (quasi-)local counterpart defined as
\begin{equation}
	\label{eq:deflocbasis}
	\tilde \frakB^\ell \coloneqq   \bigcup_{b \in \frakB}(1-\calCl)\PH b.
\end{equation}
Similarly, we also define a localized ansatz space for the practical multiscale method,
\begin{equation*}
	\tUHl \coloneqq   (1-\calCl)\PH \UH.
\end{equation*}
This choice is different from the choice in \cite{Mai21} as it employs a localization based on element corrections and additionally utilizes the operator $\PH$ in the basis definition. Henceforth, for $T \in \calT_H$, we denote by $\Lambda_{1,T}$ and  $b_{1,T}$ the constant Legendre basis function and its associated bubble function, respectively. In the following lemma, it is shown that $\PH$ only has an effect on the bubble functions {corresponding to the polynomial degree $0$, i.e., the set}
\begin{equation*}
	\frakB^0\coloneqq \{b_{1,T}\with T \in \calT_H\} \subset \frakB.
\end{equation*}
\begin{lemma}[Application of $\PH$]\label{lem:applicationPH}
	For all $T \in \calT_H$ and $j > 1$, it holds
	\begin{equation*}
		\PH b_{j,T} = b_{j,T}.
	\end{equation*}
\end{lemma}

\begin{proof}
	We consider a fixed element $T \in \calT_H$. By definition, the bubble function $b_{j,T}$ satisfies the relation $\PiH b_{j,T} = \Lambda_{j,T}$. Using this, we obtain for all $j>1$
	\begin{equation*}
		\label{eq:implementationobservation}
		\PiH^0 b_{j,T} = \tspf{b_{j,T}}{\Lambda_{1,T}}{T}\Lambda_{1,T} = \tspf{\Lambda_{j,T}}{\Lambda_{1,T}}{T}\Lambda_{1,T} = 0.
	\end{equation*}
	Hence, we have $\IH b_{j,T} = 0$, which implies $\PH b_{j,T} = b_{j,T}$ using~\cref{eq:defPH}.
\end{proof}

We show that, for any polynomial degree $p$, the construction in~\cref{eq:deflocbasis} ensures a more robust localization compared to \cite{Mai21}. Notably, for this positive effect only the basis functions corresponding to elements in $\frakB^0$ need to be modified. Their supports are enlarged by one additional layer of elements. This strategy has been used in~\cite{HauP21} to cure instabilities for the case $p = 0$; see also~\cite{AltHP21}.\medskip

The final practical multiscale method approximates the solution $u$ of~\cref{eq:prob} by performing~a Galerkin approach based on the space $\tUHl$. It seeks $\uCl \in \tUHl$ such that
\begin{equation}\label{eq:msprac}
	a(\uCl,\vCl) = (f,\vCl)_{L^2(D)}\qquad\text{for all } \vCl \in \tUHl.
\end{equation}
The next theorem proves an error estimate for the solution $\uCl$.
\begin{theorem}[A priori error estimate]\label{t:fullerror}
	Assume that $f \in H^k(\calT_H),\, k \in \N_0$ and set $s \coloneqq  \min\{k,p+1\}$. Further, let $u \in H^1_0(D)$ and $\uCl \in \tUHl$ be the solutions to \cref{eq:prob} and \cref{eq:msprac}, respectively. Then, there exists a constant $\Cel>0$, which solely depends on the regularity of $\calT_H$ and the polynomial degree $p$, such that
	\begin{align}
		\tnormf{\nabla(u-\uCl)}{D} &\leq {\Ce H^{s+1}|f|_{H^s(\calT_H)} + \Cel\ell^{d/2} \exp(-\Cd \ell)\tnormf{f}{D}}\label{eq:H1est},\\
		\tnormf{u-\uCl}{D} &\leq {\beta \big(\Ce H + \Cel\ell^{d/2} \exp(-\Cd\ell)\big)}\tnormf{\nabla(u-\uCl)}{D}\label{eq:L2est}
	\end{align}
	{with the constants $\Ce$ from \Cref{t:errODp} and $\Cd$ from \Cref{lemma:expdec}. With respect to the polynomial degree, we have $\Cel \sim (p+1)^{3}$.}
\end{theorem}

\begin{proof}
	By C\'ea's Lemma, we obtain for any $v \in \tUHl$
	\begin{equation*}
		\tnormf{\nabla (u - \uCl)}{D} \leq \sqrt{\tfrac{\beta}{\alpha}} \tnormf{\nabla (u-v)}{D}.
	\end{equation*}
	Adding and subtracting {the solution $\uC\in\tUH$ to~\cref{eq:msprob}} and employing the triangle inequality, we obtain
	\begin{equation}\label{eq:proofFinal}
		\tnormf{\nabla (u-v)}{D} \leq \tnormf{\nabla (u-\uC)}{D}+\tnormf{\nabla (\uC-v)}{D}.
	\end{equation}
	The first term is the error of the prototypical method for which we use the bound from \Cref{t:errODp}.
	{
		To bound the second term, we require the equality $\uC = (1-\calC)\PH u$, which can be shown as follows. First, we observe that
		\begin{equation*}
			(1-\calC)\PH u = (1-\calC)u - (1-\calC)(1-\PH)u = (1-\calC)u
		\end{equation*}
		since $(1-\PH)u \in \W = \ker\PiH\vert_{H^1_0(D)}$. Further, we use that any $\vC \in \tUH$ can be written as $\vC=(1-\calC)b_H$ for some $b_H \in \UH$ and calculate
		\begin{align*}
			a((1-\calC)u,\vC) &= a(u,\vC) - a(\calC u,\vC) \\&= (f,\vC)_{L^2(D)} - a(\calC u,(1-\calC)b_H)\\& = (f,\vC)_{L^2(D)},
		\end{align*}
		where we use the symmetry of $a$ and the definition of $\calC$ (cf.~\cref{eq:corr}) in the last step. Since $\uC$ is uniquely defined, we deduce that $\uC = (1-\calC)u = (1-\calC)\PH u$.
		
		We now go back to the second term in~\cref{eq:proofFinal}. Since $\PH(u - \calB_H u) = 0$ (cf.~the representation of $\PH$ in~\cref{eq:repPH}), $v \coloneqq (1-\calCl)\PH u = (1-\calCl)\PH\calB_H u \in \tUHl$ is a valid approximation of $\uC$ in the space $\tUHl$. Using this and the above equality for $\uC$, we obtain}
	\begin{align*}
		\begin{split}
			\tnormf{\nabla (\uC-v)}{D} &= \tnormf{\nabla (\calC - \calCl)\PH u}{D} \\&\leq \Cl \ell^{d/2}\exp(-\Cd \ell)\tnormf{\nabla \PH u}{D}\\
			&\leq \Cp\Cl \ell^{d/2} \exp(-\Cd\ell) \tnormf{\nabla u}{D} \\&\leq \alpha^{-1}\CF\Cp\Cl  \ell^{d/2} \exp(-\Cd\ell) \tnormf{f}{D}
		\end{split}
	\end{align*}
	employing also \Cref{lemma:locerrcorr} as well as the estimates~\cref{eq:stabu} and~\cref{eq:PHprop}. Combining the above error bounds yields the desired $H^1$-error estimate~\cref{eq:H1est}. {Moreover, the scaling of the constant $\Cel$ with respect to $p$ follows directly from the scaling of $\Cp$ and $\Cl$.}
	
	To prove \cref{eq:L2est}, we require the Aubin--Nitsche duality argument. We denote by $z_g \in H^1_0(D)$ the solution to~\cref{eq:prob} for the right-hand side $g \coloneqq u-\uCl \in L^2(D)$. Using the Galerkin orthogonality, we obtain for any $v \in \tUHl$
	\begin{align*}
		\|u-\uCl\|_{L^2(D)}^2 &= \tspf{g}{u-\uCl}{D} = a(z_g,u-\uCl) = a(z_g-v,u-\uCl).
	\end{align*}
	Choosing $v$ as the solution of the practical multiscale method~\cref{eq:msprac} for the particular right-hand side $g$ and using \cref{eq:H1est} for $s = 0$, we get
	\begin{align*}
		\|u-\uCl\|_{L^2(D)}^2 &\leq \beta \tnormf{\nabla (z_g-v)}{D}\tnormf{\nabla (u-\uCl)}{D}\\
		&\leq {\beta (\Ce H + \Cel\ell^{d/2}\exp(-\Cd \ell))}\tnormf{g}{D}\tnormf{\nabla (u-\uCl)}{D}\\
		&= {\beta (\Ce H + \Cel\ell^{d/2}\exp(-\Cd \ell))}\tnormf{u-\uCl}{D}\tnormf{\nabla (u-\uCl)}{D}.
	\end{align*}
	Dividing by $\tnormf{u-\uCl}{D}$ on both sides, the desired $L^2$-norm estimate~\cref{eq:L2est} follows.
\end{proof}

{
	\begin{remark}[Choice of $\ell$]
		The estimates in \Cref{t:fullerror} suggest to choose $\ell \gtrsim (p+1)(s+1)|\log(H)|$ in order to retain the overall convergence rates as stated in \Cref{t:errODp}. This is in line with the theoretical requirements in~\cite[Cor.~4.6]{Mai21}. We emphasize, however, that this dependence on $p$ is rather pessimistic as also observed in numerical experiments.
	\end{remark}
}

\begin{figure}[t]
	\begin{tabularx}{.9\textwidth}{@{}YY@{}}
		$p = 0,1$ & $p = 2,3$\\
		\includegraphics[trim=0 0 50 20, clip,width=\linewidth]{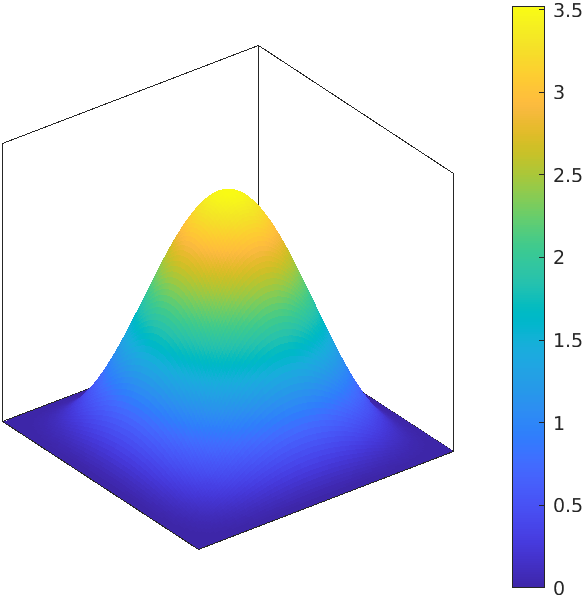}&
		\includegraphics[trim=0 0 50 20, clip,width=\linewidth]{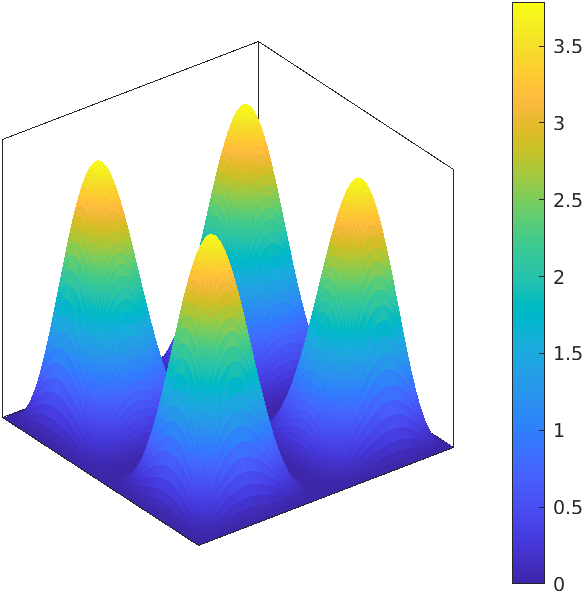}
	\end{tabularx}
	\hfill
	\begin{tabularx}{.065\textwidth}{@{}Y@{}}
		\includegraphics[trim=240 0 0 0, clip,width=\linewidth]{figures/bubble_p0}
	\end{tabularx}
	\caption{Bubble functions $b_{1,T}$ for polynomial degrees $p = 0,1$ (left) and polynomial degrees $p = 2,3$ (right).}
	\label{fig:bubbles}
\end{figure}

\section{Numerical examples}\label{s:numerics}

This section numerically investigates the proposed improved high-order method. For all numerical experiments, we use the domain $D = (0,1)^2$ and consider a coarse Cartesian mesh $\calT_H$ of $D$. Note that, henceforth, the mesh size denotes the element side-length instead of the diameter. The subsequent numerical experiments require an explicit construction of bubble functions. One possible choice is outlined in the following remark.
{
	\begin{remark}[Construction of bubbles]\label{remark:constrbubbles}
		For the construction of the bubble functions associated to an element $T \in \calT_H$, we introduce the function $\theta \in H^1_0(T)$ obtained by multiplying all nodal $\mathcal Q_1$ basis functions within $T$. The bubble functions $b_{j,T}$ are then obtained by the ansatz $b_{j,T} = \sum_{i = 1}^{N} c_i\, \theta\Lambda_{i,T}$, where the coefficients are chosen such that
		\begin{equation*}
			\PiH b_{j,T} = \Lambda_{j,T}.
		\end{equation*}
		This condition translates to a $N \times N$ linear system, which is invertible under the condition that $\{\PiH (\theta\Lambda_{j,T})\with j = 1,\dots,N\}$ are linearly independent. The linear independence can be proved utilizing the positivity of $\theta$ in the interior of $T$. 
		
		We emphasize that, by \cref{lem:applicationPH} and the localized counterpart  of \cref{eq:saddlepoint}, only the zero-order bubble functions are explicitly required for the computation of the basis functions~\cref{eq:deflocbasis}. The remaining basis functions can be computed as in~\cite{Mai21}. Zero-order bubble functions corresponding to different polynomial degrees which fulfill the above conditions are depicted in \cref{fig:bubbles}. We note that a symmetry argument shows that the bubble functions for consecutive polynomial degrees ($p = k$ and $p = k+1$, $k = 0,2,4,\dots $) coincide.
	\end{remark}
}
{To obtain a practically computable method, we also} need to discretize the infinite-dimensional localized corrections defined in~\cref{eq:locelemcorr}. For this, we employ the $\mathcal Q_1$-finite element method on a fine mesh of the respective patches with mesh size $h = 2^{-9}$ obtained by uniform refinement of the coarse patch mesh. For details regarding the discretization of these problems, we refer to~\cite[Sec.~4.3]{Mai21}. We emphasize that also the construction of the bubbles in \cref{remark:constrbubbles} carries over to the fully discrete setting. All errors are computed with respect to the $\mathcal Q_1$-finite element reference solution on the global Cartesian mesh with mesh size $h = 2^{-9}$. All considered errors are relative errors measured in the energy norm
\begin{equation*}
	\|\cdot \|_a \coloneqq \tnormf{A^{1/2}\nabla \cdot}{D}.
\end{equation*}

\begin{figure}[t]
	\begin{tabularx}{\textwidth}{@{}YY@{}}
		\includegraphics[width=\linewidth]{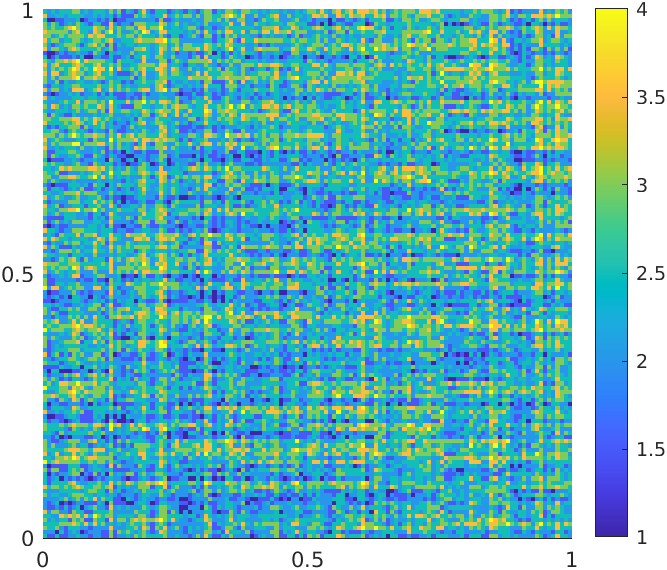} & 	\includegraphics[width=\linewidth]{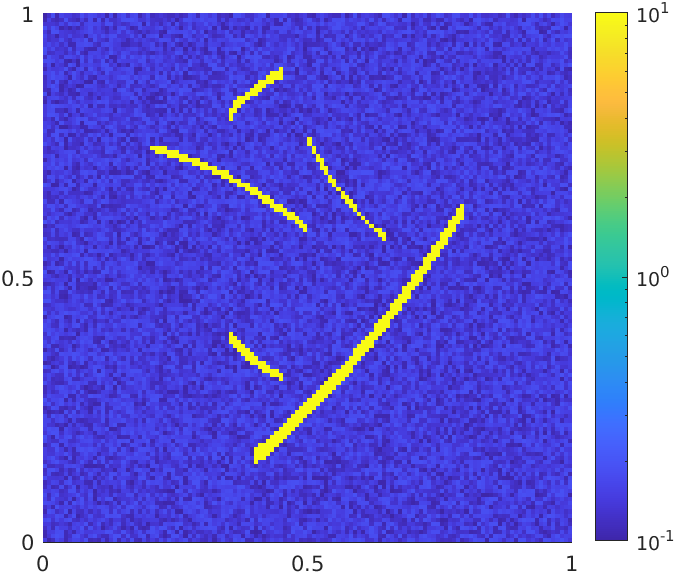}
	\end{tabularx}
	\caption{Coefficients $A_1$, $A_2$ used in the numerical experiments.}
	\label{fig:coeff}
\end{figure}

For the numerical experiments, we use the two scalar coefficients $A_1$, $A_2$ as depicted in \Cref{fig:coeff}. Both coefficients are piecewise constant with respect to the Cartesian mesh~$\calT_{2^{-7}}$. The coefficient $A_1$ is the sum of random contribution on several scales with~$2^{-7}$ being the smallest scale. It varies within the interval $[1,4]$. The coefficient $A_2$ features some channels of conductivity 10 which are placed on top of a noisy background. On each element of~$\calT_{2^{-7}}$, the noisy background takes values generated by realizations of independent random variables which are uniformly distributed in the interval $[0.1,1]$.
As right-hand side, we consider the smooth function
\begin{equation*}
	f(x_1,x_2) = (x_1+\cos(3\pi x_1))\cdot x_2^3.
\end{equation*}

\begin{figure}[h]
	\begin{tabularx}{\textwidth}{@{}YY@{}}
		$p = 0$ & $p = 1$\\
		\includegraphics[width=\linewidth]{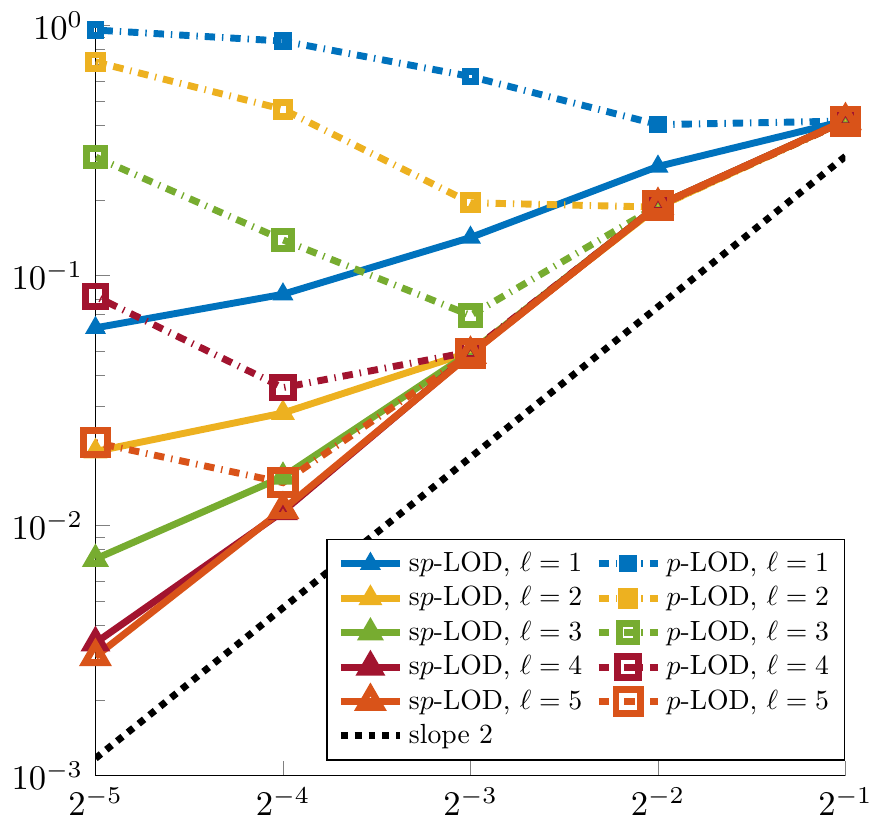}&
		\includegraphics[width=\linewidth]{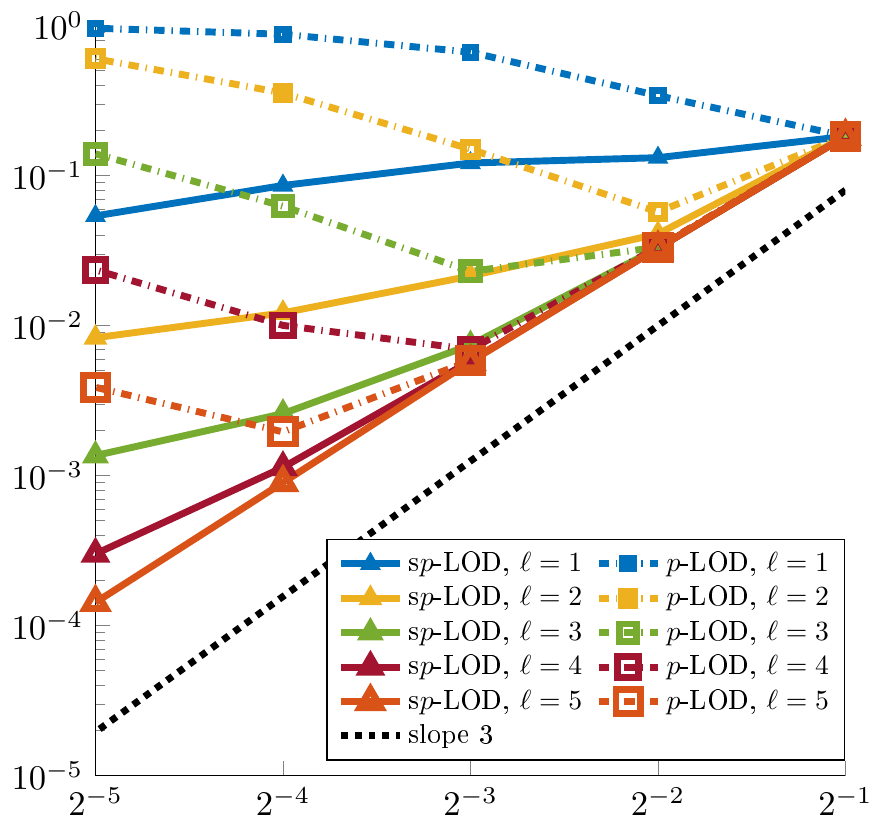}\\[2ex]
		$p = 2$ & $p = 3$\\
		\includegraphics[width=\linewidth]{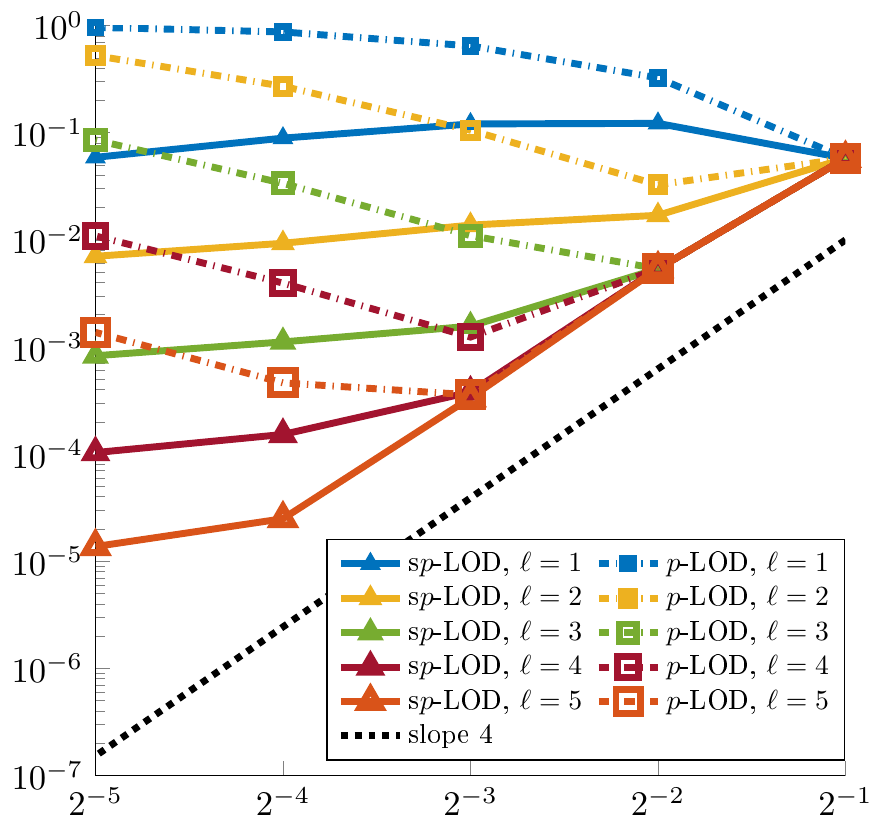}&
		\includegraphics[width=\linewidth]{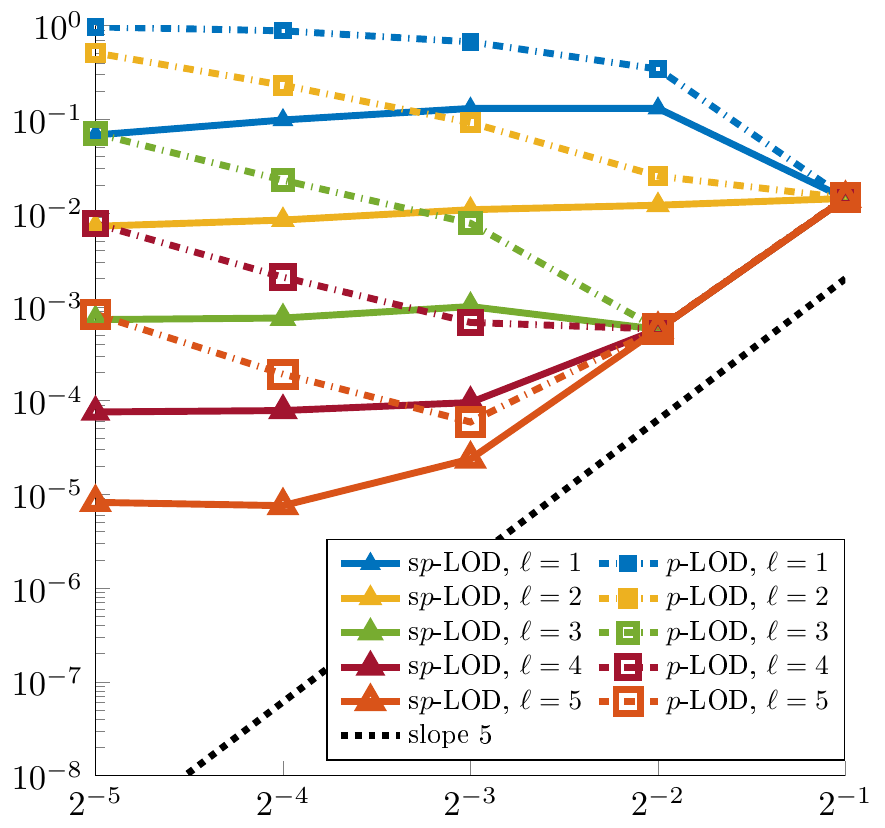}
	\end{tabularx}
	\caption{Relative errors with respect to $H$ for the coefficient $A_1$, different values of $\ell$, and polynomial degrees $p = 0,\dots,3$ (top left to bottom right).}
	\label{fig:convA1}
\end{figure}

For the first numerical experiment, we consider the diffusion coefficient $A_1$ and investigate the convergence of the proposed stabilized high-order method (referred to as s$p$-LOD) for polynomial degrees $p= 0,\dots,3$. The results are compared to the localized high-order method in~\cite{Mai21} (referred to as $p$-LOD). For fixed oversampling parameters, \Cref{fig:convA1} plots the relative energy errors of both methods as a function of the coarse mesh size $H$. As reference, also a line indicating the expected rate of convergence is shown.
Provided that the oversampling parameter $\ell$ is chosen large enough, one observes convergence of order $p+2$ for the methods of degree $p$. This numerically confirms the theoretical prediction from \Cref{t:errODp}. Most importantly, \Cref{fig:convA1} confirms the improved stability properties of the s$p$-LOD compared to the $p$-LOD. For fixed oversampling parameters, the error of the s$p$-LOD first decreases and then stagnates as $H$ is decreased. In contrast, one observes an increasing error of the $p$-LOD for mesh sizes below a certain critical threshold depending on the oversampling parameter, which is in line with the theoretical results in~\cite[Thm.~4.4]{Mai21}. Note that for both localization strategies, the stagnation (respectively deterioration) of the error can be avoided if the oversampling parameter is appropriately increased when $H$ is decreased (i.e., logarithmically in $1/H$; cf.~also \Cref{t:fullerror}).

{Next, we specifically investigate the localization errors. Therefore, we choose the right-hand side $f \equiv 1$ for which the first term on the right-hand side of the error estimate~\cref{eq:H1est} vanishes and, hence, only the localization error remains. For different choices of $\ell$ and $p$, we numerically investigate the decay of the localization error.
	\begin{figure}[h]
		\begin{tabularx}{\textwidth}{@{}YY@{}}
			\includegraphics[width=\linewidth]{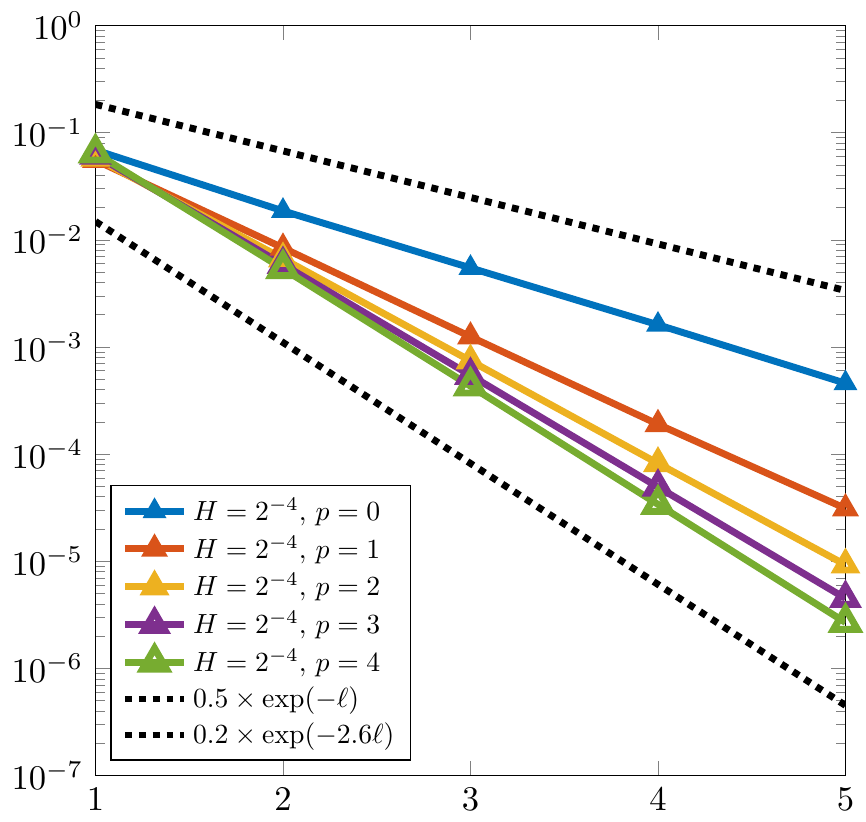}&
			\includegraphics[width=\linewidth]{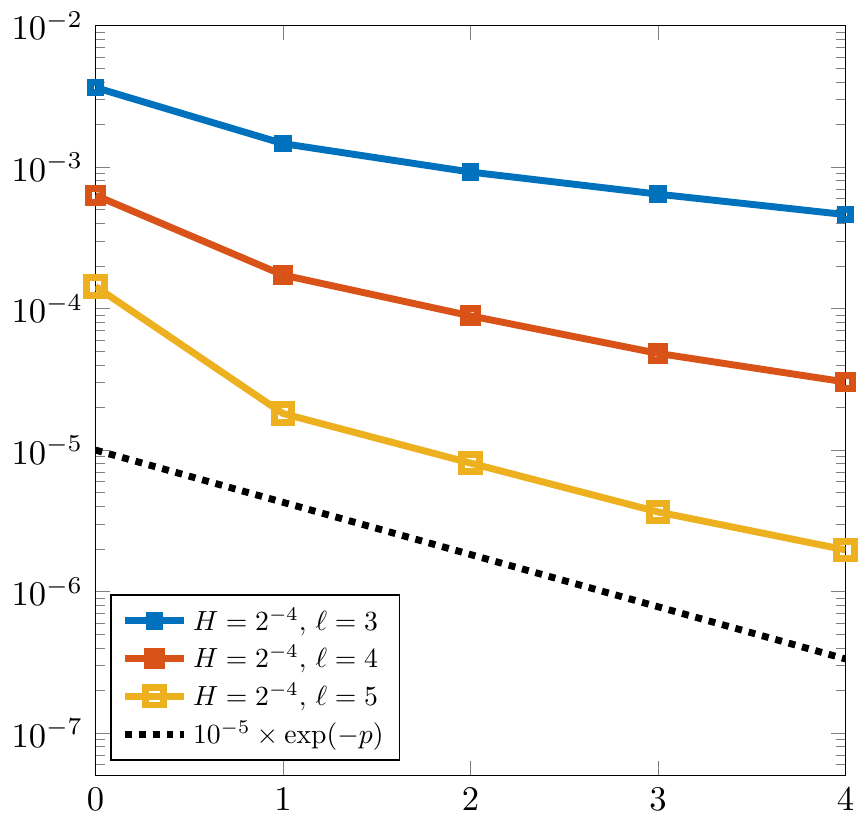}
		\end{tabularx}
		\caption{{Relative localization errors for the coefficient $A_1$ and $H = 2^{-4}$ for fixed $p$ and varying $\ell$ (left) and for fixed $\ell$ and varying $p$ (right).}}
		\label{fig:dec}
	\end{figure}
	In \cref{fig:dec}, one observes that the localization error exponentially decays as $\ell$ is increased which numerically confirms~\cref{t:fullerror}. Further, one observes that for larger values of~$p$, the decay behavior is further improved. This can also be seen in \cref{fig:dec}, where the errors are plotted with respect to $p$ for different fixed values of~$\ell$. This is in line with the numerical results in~\cite{Mai21} and indicates that the scaling of $\Cd$ in~\Cref{lemma:expdec} with respect to $p$ is most likely not sharp.}

\begin{figure}[h]
	\begin{tabularx}{.92\textwidth}{@{}YY@{}}
		$p = 0$ & $p = 1$\\
		\includegraphics[trim=0 0 50 20, clip,width=\linewidth]{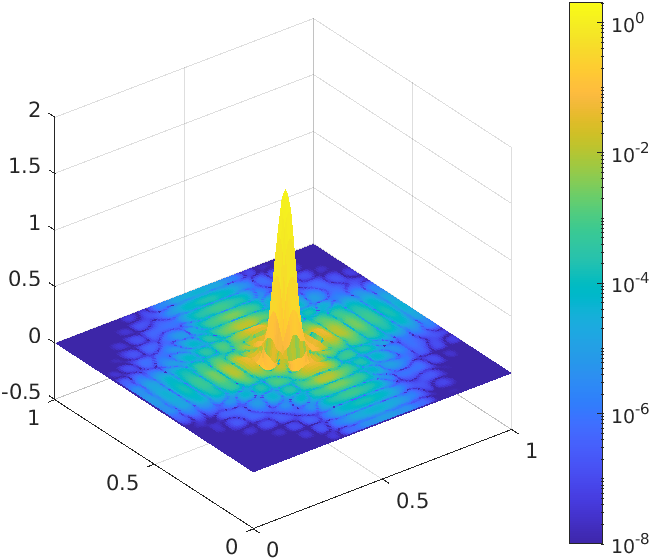}&
		\includegraphics[trim=0 0 50 20, clip,width=\linewidth]{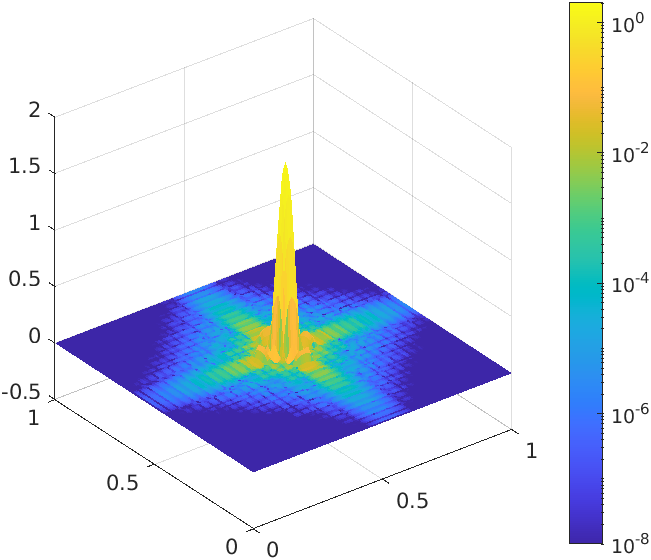}\\[2ex]
		$p = 2$ & $p = 3$\\
		\includegraphics[trim=0 0 50 20, clip,width=\linewidth]{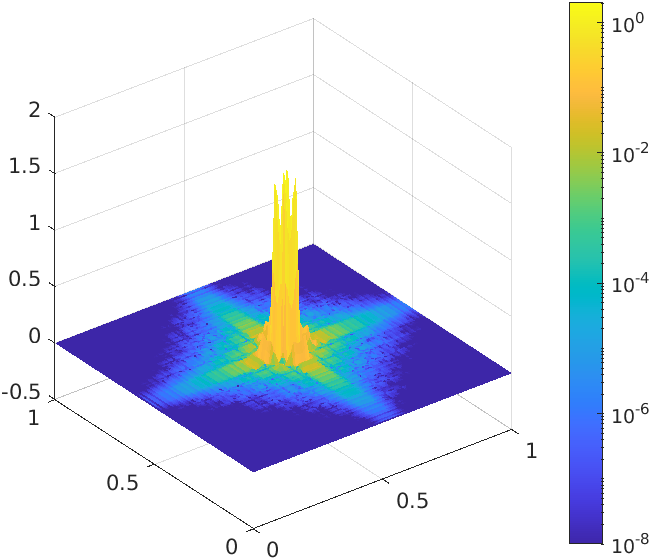}&
		\includegraphics[trim=0 0 50 20, clip,width=\linewidth]{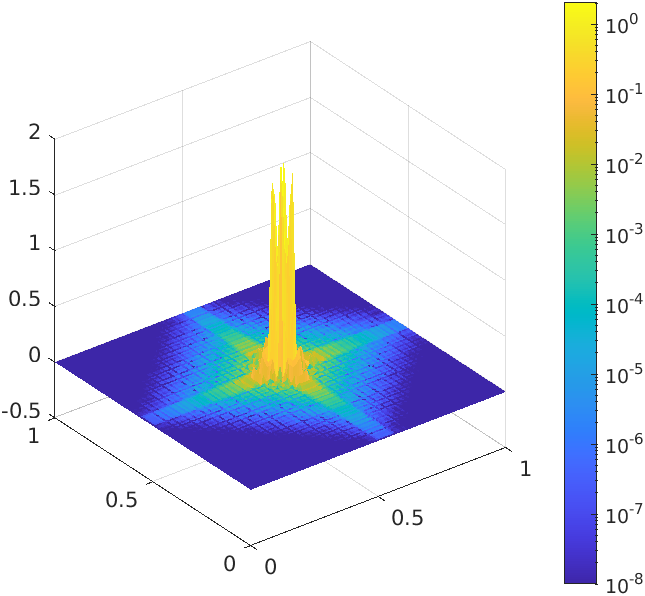}
	\end{tabularx}
	\begin{tabularx}{.065\textwidth}{@{}Y@{}}
		\includegraphics[trim=272 0 0 0, clip,width=\linewidth]{figures/basisfunp0revision}
	\end{tabularx}
	\caption{Ideal basis functions $(1-\calC)b_{1,T}$ for the polynomial degrees $p = 0,\dots,3$ (top left to bottom right).}
	\label{fig:basisfun}
\end{figure}

Next, we depict selected prototypical basis functions. Therefore, we choose a fixed element $T$ of the coarse Cartesian mesh $\calT_{2^{-4}}$. For polynomial degrees $p = 0,\dots,3$, \Cref{fig:basisfun} shows plots of the discretized counterparts of the prototypical basis functions $(1-\calC)b_{1,T}$, where $b_{1,T} \in H^1_0(T)$ is the bubble function satisfying $\PiH b_{1,T} = \Lambda_{1,T}$ with $\Lambda_{1,T}$ denoting the $L^2$-normalized characteristic function of the element $T$ which is the first Legendre basis function. Note that the $p$-LOD and the s$p$-LOD have the same prototypical basis functions as both methods differ only with respect to the localization strategy. For all plotted basis functions in \Cref{fig:basisfun}, we use the same color coding indicating the absolute value of the respective basis functions. We employ a logarithmic color map which we restrict to the interval $[10^{-8},3.5]$. This means that all values below $10^{-8}$ are displayed in dark blue indiscriminate of their actual value.
In \Cref{fig:basisfun}, one observes that, for larger polynomial degrees $p$, the basis functions tend to be more localized, i.e., a more rapid decay can be  observed. This is in line with the results in~\cite{Mai21}, where improved localization is observed for larger polynomial degrees.

\begin{figure}[h]
	\centering
	\includegraphics[width=.475\textwidth]{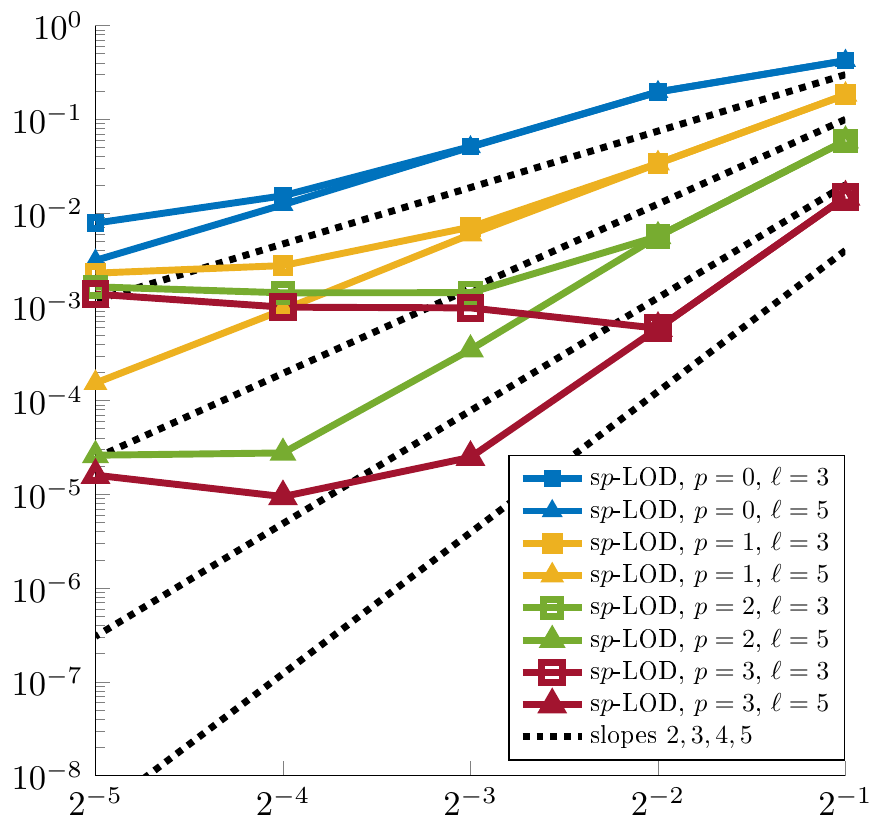}
	\caption{Convergence plot for  the coefficient $A_2$ and selected values of~$\ell$ and $p$. }
	\label{fig:convA2}
\end{figure}

In the second numerical experiment, we investigate the convergence of the s$p$-LOD for the more challenging diffusion coefficient $A_2$. For polynomial degrees $p = 0,\dots,3$,
\Cref{fig:convA2} plots the errors for oversampling parameters $\ell = 3,5$ as a function of~$H$. As reference, lines indicating the respective expected convergence rates are depicted.
One observes that the errors plotted in \Cref{fig:convA2} are qualitatively similar to the ones in \Cref{fig:convA1}. One obtains the expected optimal convergence rates of order $p+2$ provided that the oversampling parameter is chosen sufficiently large. Again, the method demonstrates its improved stability properties, as the error first decreases and then stagnates for decreasing~$H$ and does not show a deteriorating behavior. As mentioned above, the convergence rate of order $p+2$ can be preserved if the oversampling parameter $\ell$ is logarithmically coupled to the mesh size~$H$ as predicted by the theoretical result in \Cref{t:fullerror}.

\section{Conclusions}\label{s:conclusions}

In this work, we have presented an improved version of the high-order localized multiscale method presented in~\cite{Mai21} which can be used to numerically solve elliptic multiscale problems with rough coefficients. The method only requires minimal assumptions on the coefficient and achieves higher-order approximation rates exploiting regularity properties of the right-hand side only. The method is defined in the spirit of the Localized Orthogonal Decomposition methodology and is built upon appropriate high-order bubble functions. An improved localization strategy has been presented that cures the issue of deteriorating errors for smaller mesh sizes as theoretically and practically observed in~\cite{Mai21}.
It slightly enlarges the support of the lowest-order basis functions and thereby notably achieves a stabilization for arbitrary polynomial degree.
We have provided a complete a priori error analysis and have demonstrated the performance and the improvements of the method in a set of numerical experiments.  

\bibliographystyle{siamplain}
\bibliography{bib}

\begin{thebibliography}{AHPV13}

\bibitem[AB05]{AllB05}
G.~Allaire and R.~Brizzi.
\newblock A multiscale finite element method for numerical homogenization.
\newblock {\em Multiscale Model. Simul.}, 4(3):790--812, 2005.

\bibitem[AB12]{AbdB12}
A.~Abdulle and Y.~Bai.
\newblock Reduced basis finite element heterogeneous multiscale method for
  high-order discretizations of elliptic homogenization problems.
\newblock {\em J. Comput. Phys.}, 231(21):7014--7036, 2012.

\bibitem[AHP21]{AltHP21}
R.~Altmann, P.~Henning, and D.~Peterseim.
\newblock Numerical homogenization beyond scale separation.
\newblock {\em Acta Numer.}, 30:1--86, 2021.

\bibitem[AHPV13]{AraHPV13}
R.~Araya, C.~Harder, D.~Paredes, and F.~Valentin.
\newblock Multiscale hybrid-mixed method.
\newblock {\em SIAM J. Numer. Anal.}, 51(6):3505--3531, 2013.

\bibitem[BCO94]{BabCO94}
I.~Babu\v{s}ka, G.~Caloz, and J.~E. Osborn.
\newblock Special finite element methods for a class of second order elliptic
  problems with rough coefficients.
\newblock {\em SIAM J. Numer. Anal.}, 31(4):945--981, 1994.

\bibitem[BL11]{BabL11}
I.~Babu\v{s}ka and R.~Lipton.
\newblock Optimal local approximation spaces for generalized finite element
  methods with application to multiscale problems.
\newblock {\em Multiscale Model. Simul.}, 9(1):373--406, 2011.

\bibitem[BO83]{BabO83}
I.~Babu\v{s}ka and J.~E. Osborn.
\newblock Generalized finite element methods: their performance and their
  relation to mixed methods.
\newblock {\em SIAM J. Numer. Anal.}, 20(3):510--536, 1983.

\bibitem[BO00]{BabO00}
I.~Babu\v{s}ka and J.~E. Osborn.
\newblock Can a finite element method perform arbitrarily badly?
\newblock {\em Math. Comp.}, 69(230):443--462, 2000.

\bibitem[Bre94]{Bre94}
S.~C. Brenner.
\newblock Two-level additive {S}chwarz preconditioners for nonconforming finite
  elements.
\newblock {\em Contemp. Math.}, 180:9--14, 1994.

\bibitem[CEL19]{CicEL19}
M.~Cicuttin, A.~Ern, and S.~Lemaire.
\newblock A hybrid high-order method for highly oscillatory elliptic problems.
\newblock {\em Comput. Methods Appl. Math.}, 19(4):723--748, 2019.

\bibitem[EG17]{ErnG17}
A.~Ern and J.-L. Guermond.
\newblock Finite element quasi-interpolation and best approximation.
\newblock {\em ESAIM Math. Model. Numer. Anal.}, 51(4):1367--1385, 2017.

\bibitem[EGH13]{EfeGH13}
Y.~Efendiev, J.~Galvis, and T.~Y. Hou.
\newblock Generalized multiscale finite element methods ({GM}s{FEM}).
\newblock {\em J. Comput. Phys.}, 251:116--135, 2013.

\bibitem[GGS12]{GraGS12}
L.~Grasedyck, I.~Greff, and S.~Sauter.
\newblock The {AL} basis for the solution of elliptic problems in heterogeneous
  media.
\newblock {\em Multiscale Model. Simul.}, 10(1):245--258, 2012.

\bibitem[HP13]{HenP13}
P.~Henning and D.~Peterseim.
\newblock Oversampling for the multiscale finite element method.
\newblock {\em Multiscale Model. Simul.}, 11(4):1149--1175, 2013.

\bibitem[HP22]{HauP21}
M.~Hauck and D.~Peterseim.
\newblock Multi-resolution localized orthogonal decomposition for {H}elmholtz
  problems.
\newblock {\em Multiscale Model. Simul.}, 20(2):657--684, 2022.

\bibitem[HPV13]{HarPV13}
C.~Harder, D.~Paredes, and F.~Valentin.
\newblock A family of multiscale hybrid-mixed finite element methods for the
  {D}arcy equation with rough coefficients.
\newblock {\em J. Comput. Phys.}, 245:107--130, 2013.

\bibitem[HSS02]{HouSS02}
P.~Houston, C.~Schwab, and E.~S\"{u}li.
\newblock Discontinuous {$hp$}-finite element methods for
  advection-diffusion-reaction problems.
\newblock {\em SIAM J. Numer. Anal.}, 39(6):2133--2163, 2002.

\bibitem[HZZ14]{HesZZ14}
J.~S. Hesthaven, S.~Zhang, and X.~Zhu.
\newblock High-order multiscale finite element method for elliptic problems.
\newblock {\em Multiscale Model. Simul.}, 12(2):650--666, 2014.

\bibitem[LMT12]{LiPT12}
R.~Li, P.~Ming, and F.~Tang.
\newblock An efficient high order heterogeneous multiscale method for elliptic
  problems.
\newblock {\em Multiscale Model. Simul.}, 10(1):259--283, 2012.

\bibitem[Mai20]{Mai20}
R.~Maier.
\newblock {\em Computational Multiscale Methods in Unstructured Heterogeneous
  Media}.
\newblock PhD thesis, University of Augsburg, 2020.

\bibitem[Mai21]{Mai21}
R.~Maier.
\newblock A high-order approach to elliptic multiscale problems with general
  unstructured coefficients.
\newblock {\em SIAM J. Numer. Anal.}, 59(2):1067--1089, 2021.

\bibitem[MP14]{MalP14}
A.~M{\aa}lqvist and D.~Peterseim.
\newblock Localization of elliptic multiscale problems.
\newblock {\em Math. Comp.}, 83(290):2583--2603, 2014.

\bibitem[MP20]{MalP20}
A.~M\aa{}lqvist and D.~Peterseim.
\newblock {\em Numerical homogenization by localized orthogonal decomposition},
  volume~5 of {\em SIAM Spotlights}.
\newblock Society for Industrial and Applied Mathematics (SIAM), Philadelphia,
  PA, 2020.

\bibitem[OS19]{OwhS19}
H.~Owhadi and C.~Scovel.
\newblock {\em Operator-adapted wavelets, fast solvers, and numerical
  homogenization}, volume~35 of {\em Cambridge Monographs on Applied and
  Computational Mathematics}.
\newblock Cambridge University Press, Cambridge, 2019.

\bibitem[Osw93]{Osw93}
P.~Oswald.
\newblock On a {BPX}-preconditioner for {P}1 elements.
\newblock {\em Computing}, 51(2):125--133, 1993.

\bibitem[Owh17]{Owh17}
H.~Owhadi.
\newblock Multigrid with rough coefficients and multiresolution operator
  decomposition from hierarchical information games.
\newblock {\em SIAM Rev.}, 59(1):99--149, 2017.

\bibitem[Wey16]{Wey16}
M.~Weymuth.
\newblock {\em Adaptive local basis for elliptic problems with
  {$L^\infty$}-coefficients}.
\newblock PhD thesis, University of Zurich, 2016.

\end{thebibliography}

\appendix
\section{Proofs of \cref{lemma:expdec,lemma:locerrcorr}}
\label{sec:appendix}

\begin{proof}[Proof of \Cref{lemma:expdec}]
	This proof uses classical LOD techniques and is similar to the proofs of \cite[Thm.~4.1]{MalP20} and \cite[Lem.~5.2]{HauP21}.
	For shorter notation, let us denote $\varphi\coloneqq \calCT v$. We define the finite element cut-off function $\eta\in W^{1,\infty}(D,[0,1])$ such that
	\begin{alignat*}{2}
		&\eta \equiv 0 \quad &&\text{ in }\Nb^{\ell-1}(T),\\
		&\eta \equiv 1 \quad &&\text{ in }D\backslash \Nb^{\ell}(T),\\
		&0 \leq \eta \leq 1 \qquad &&\text{ in } R\coloneqq \Nb^{\ell}(T)\backslash \Nb^{\ell-1}(T).
	\end{alignat*}
	with $\|\nabla \eta\|_{L^\infty(D)}\leq C_\eta H^{-1}$ for some $C_\eta>0$. We estimate
	\begin{align}
		\alpha\tnorm{\nabla \varphi}{D\backslash \Nb^\ell(T)}^2 &\leq  (A\nabla\varphi,\eta\nabla\varphi)= a(\varphi,\eta\varphi) - (A\nabla \varphi,\varphi \nabla \eta)\notag\\
		&\leq \big\vert a(\varphi,(1-\calB_H)(\eta\varphi))\big\vert   +\big\vert a(\varphi,\calB_H(\eta\varphi))\big\vert + \big\vert (A\nabla \varphi,\varphi \nabla \eta)\big\vert\label{eq:refdec}
	\end{align}
	For the first term we utilize definition \cref{eq:elemcorrdef} and the fact that $w\coloneqq (1-\calB_H)(\eta\varphi) \in \W$ and $w\equiv 0$ in $T$. For $\ell\geq 1$, this yields
	\begin{equation*}
		a(\varphi,(1-\calB_H)(\eta\varphi)) = a_T(\varphi,w) = 0.
	\end{equation*}
	Using that $\calB_H(\eta\varphi) \equiv 0$ in $D\backslash \Nb^{\ell}(T)$ and employing the estimates \cref{eq:inv} and \cref{eq:approx}, we obtain for the second term in \cref{eq:refdec}
	\begin{align*}
		\big\vert a(\varphi,\calB_H(\eta\varphi))\big\vert&\leq \Cbo \beta H^{-1}\tnorm{\nabla\varphi}{R}\tnorm{\eta\varphi}{R}\\
		&\leq \Cbo \beta H^{-1}\tnorm{\nabla\varphi}{R}\tnorm{(1-\PiH)\varphi}{R}\leq \Cbo\Ca\beta \tnorm{\nabla\varphi}{R}^2.
	\end{align*}
	For the third term in \cref{eq:refdec}, we obtain with similar arguments
	\begin{align*}
		\big\vert (A\nabla \varphi,\varphi \nabla \eta)\big\vert \leq C_\eta \beta  H^{-1}\tnorm{\nabla \varphi}{R}\tnorm{(1-\PiH)\varphi}{R}\leq C_\eta \Ca\beta \tnorm{\nabla\varphi}{R}^2.
	\end{align*}
	Combining the previous estimates yields
	\begin{align*}
		\tnorm{\nabla\varphi}{D\backslash \Nb^\ell(T)}^2\leq C \tnorm{\nabla\varphi}{R}^2 &=   C \tnorm{\nabla\varphi}{D\backslash \Nb^{\ell-1}(T)} - C \tnorm{\nabla\varphi}{D\backslash \Nb^\ell(T)}^2\\
		\Leftrightarrow\; \big(1+C\big)\tnorm{\nabla\varphi}{D\backslash \Nb^\ell(T)}^2&\leq C \tnorm{\nabla\varphi}{D\backslash \Nb^{\ell-1}(T)}^2\\
		\Leftrightarrow \; \tnorm{\nabla\varphi}{D\backslash \Nb^\ell(T)} &\leq \sqrt{\frac{C}{1+C}}\, \tnorm{\nabla\varphi}{D\backslash \Nb^{\ell-1}(T)},
	\end{align*}
	where $C = \Ca(\Cbo+C_\eta)\beta\alpha^{-1}$. Defining $\Cd \coloneqq \frac12\big|\log\big(\tfrac{C}{1+C}\big)\big|>0$ and iterating the last inequality, we finally obtain
	\begin{equation*}
		\tnorm{\nabla\varphi}{D\backslash \Nb^{\ell}(T)}\leq \exp(-\Cd \ell )\tnormf{\nabla \varphi}{D}.
	\end{equation*}
	{Regarding the scaling of $\Cd$ with respect to $p$, we point out that the scaling of $C$ above is mainly determined by the factor $\Ca\Cbo$ and, thus, $\big|\log\big(\tfrac{C}{1+C}\big)\big|$ scales at most like ${(p+1)^{-1}}$.}
\end{proof}

For the proof of \Cref{lemma:locerrcorr}, we infer the following intermediate result.
\begin{lemma}[Localization error of element corrections]\label{lemma:locerrelemcorr}
	There exists $\Ct>0$ independent of $H,\ell,T$ such that for all $T\in\calT_H$, $v\in H^1_0(D)$, and $\ell \in\mathbb{N}$
	\begin{equation*}
		\tnormf{\nabla(\calCTl v-\calCT v)}{D} \leq \Ct \exp(-\Cd\ell) \tnormf{\nabla \calCT v}{D},
	\end{equation*}
	where $\Cd$ is the constant from  \Cref{lemma:expdec} {and $\Ct \leq C (p+1)$}.
\end{lemma}
\begin{proof}
	This proof is again similar to the proofs of \cite[Cor.~4.2]{MalP20} and  \cite[Lem.~A.1]{HauP21}.
	We abbreviate $\varphi\coloneqq \calCT v$, $\varphi^\ell\coloneqq \calCTl v$ and consider an arbitrary $w\in \W_T^\ell$. Using the Galerkin orthogonality, we obtain
	\begin{align*}
		\alpha \tnormf{\nabla (\varphi^\ell - \varphi)}{D}^2&\leq a(\varphi^\ell - \varphi,\varphi^\ell - \varphi) = a(\varphi^\ell - \varphi,w - \varphi)\\
		&\leq \beta \tnormf{\nabla (\varphi^\ell - \varphi)}{D}\tnormf{\nabla(w - \varphi)}{D}.
	\end{align*}
	We define a particular function $w \coloneqq  (1-\calB_H)(\eta\varphi)$ with
	the  finite element cut-off function $\eta\in W^{1,\infty}(D,[0,1])$ such that
	\begin{alignat*}{2}
		&\eta \equiv 1 \quad &&\text{ in }\Nb^{\ell-1}(T),\\
		&\eta \equiv 0 \quad &&\text{ in }D\backslash \Nb^{\ell}(T),\\
		&0 \leq \eta \leq 1 \qquad &&\text{ in } R\coloneqq  \Nb^{\ell}(T)\backslash \Nb^{\ell-1}(T)
	\end{alignat*}
	and $\|\nabla \eta\|_{L^\infty(D)}\leq C_\eta H^{-1}$. Since $\varphi \in \W$ and $w \in \W_T^\ell$, we obtain using the estimates \cref{eq:inv} and \cref{eq:approx}
	\begin{align*}
		\tnormf{\nabla (\varphi^\ell - \varphi)}{D}&\leq \beta \alpha^{-1} \tnorm{\nabla (w - \varphi)}{D} \\
		&= \beta\alpha^{-1} \tnormf{\nabla ((1-\calB_H)(1-\eta)\varphi)}{D\backslash \Nb^{\ell-1}(T)}\\
		&\leq (1+\Ca (C_\eta + \Cbo))\beta \alpha^{-1}\tnorm{\nabla \varphi}{D\backslash \Nb^{\ell-1}(T)}\\
		&\leq \Ct \exp(-\Cd \ell) \tnorm{\nabla\varphi}{D},
	\end{align*}
	where $\Ct \coloneqq (1+\Ca (C_\eta + \Cbo))\beta \alpha^{-1}\exp(\Cd)$. In the previous estimate, we have used that
	\begin{align}
		\begin{split}
			&\tnormf{\nabla ((1-\calB_H)(1-\eta)\varphi)}{D\backslash \Nb^{\ell-1}(T)}\\
			&\qquad \leq \tnorm{\nabla ((1-\eta)\varphi)}{D\backslash \Nb^{\ell-1}(T)} + \Cbo H^{-1}\tnorm{ (1-\eta)\varphi}{D\backslash \Nb^{\ell-1}(T)}\\
			&\qquad \leq \tnorm{\nabla \varphi}{D\backslash \Nb^{\ell-1}(T)} +(C_\eta + \Cbo)H^{-1}\tnorm{\varphi}{D\backslash \Nb^{\ell-1}(T)}\\
			&\qquad \leq(1+\Ca (C_\eta + \Cbo) )\tnorm{\nabla \varphi}{D\backslash \Nb^{\ell-1}(T)},
		\end{split}
		\label{eq:ref}
	\end{align}
	using once again the estimates \cref{eq:inv} and \cref{eq:approx}. {As in the proof of \Cref{lemma:expdec}, we observe that the scaling of $\Ct$ with respect to $p$ is determined by the scaling of $\Ca\Cbo$, which proves the stated scaling of $\Ct$.}
\end{proof}

\begin{proof}[Proof of \Cref{lemma:locerrcorr}]
	This proof is similar to the proof of \cite[Thm.~4.3]{MalP20} and \cite[Lem.~5.4]{HauP21}. We denote $z\coloneqq \calCl v - \calC v$ and, for $T\in\calT_H$$, z_T \coloneqq \calCTl v- \calCT v$. It holds
	\begin{equation*}
		\alpha \tnorm{\nabla z}{D}^2\leq a(z,z) = \sum_{T\in \calT_H} a(z_T,z).
	\end{equation*}
	Henceforth, we fix an element $T\in \calT_H$ and define the finite element the cut-off function \linebreak $\eta\in W^{1,\infty}(D,[0,1])$ such that
	\begin{alignat*}{2}
		&\eta \equiv 0 \quad &&\text{ in }\Nb^{\ell}(T),\\*
		&\eta \equiv 1 \quad &&\text{ in }D\backslash \Nb^{\ell+1}(T),\\*
		&0 \leq \eta \leq 1 \qquad &&\text{ in } R\coloneqq \Nb^{\ell+1}(T)\backslash \Nb^{\ell}(T)
	\end{alignat*}
	and $\|\nabla \eta\|_{L^\infty(D)}\leq C_\eta H^{-1}$. Using that  $\supp((1-\calB_H)(\eta z))\subset D\backslash \Nb^\ell(T)$ and $(1-\calB_H)(\eta z)\in \W$ and employing definition \cref{eq:elemcorrdef}, we get
	\begin{equation*}
		a(z_T,(1-\calB_H)(\eta z)) = -a(\calCT v,(1-\calB_H)(\eta z)) = 0.
	\end{equation*}
	Hence, using that $z = (1-\calB_H)z$ since $z\in \W$, we arrive at
	\begin{equation*}
		a(z_T,z) = a(z_T,z-(1-\calB_H)(\eta z)) = a(z_T,(1-\calB_H)((1-\eta)z)).
	\end{equation*}
	Proceeding similarly as in estimate \cref{eq:ref}, we obtain
	\begin{equation*}
		\big\vert a(z_T,(1-\calB_H)((1-\eta)z))\big\vert \leq C \beta \tnormf{\nabla z_T}{D}\tnormf{\nabla z}{\Nb^{\ell+1}(T)}.
	\end{equation*}
	with $C = (1+\Ca (C_\eta + \Cbo) )$.
	Note that the element corrections satisfy the estimate $\tnormf{\nabla \calCT v}{D}\leq  \beta \alpha^{-1} \tnormf{\nabla v}{T}$, which follows from
	\begin{align*}
		\alpha \tnorm{\nabla \calCT v}{D}^2\leq a(\calCT v,\calCT v) = a_T(v,\calCT v)
		\leq \beta \tnormf{\nabla v}{T}\tnormf{\nabla \calCT v}{D}
	\end{align*}
	using the definition \cref{eq:elemcorrdef}.
	Making use of the above inequalities, \Cref{lemma:locerrelemcorr}, the discrete Cauchy--Schwarz inequality and the finite finite overlap of the patches, we obtain after summation over all elements
	\begin{align*}
		\alpha \tnorm{\nabla z}{D}^2&\leq \sum_{T\in \calT_H} a(z_T,z)\leq C \beta  \sum_{T\in \calT_H} \tnormf{\nabla z_T}{D}\tnormf{\nabla z}{\Nb^{\ell+1}(T)}\\
		&\leq C  \Ct \beta  \exp(-\Cd \ell) \sum_{T\in \calT_H}\tnorm{\nabla \calCT v}{D}\vnorm{z}{\Nb^{\ell+1}(T)}\\
		&\leq C \Ct \beta^2 \alpha^{-1} \exp(-\Cd \ell)\sum_{T\in \calT_H} \tnormf{\nabla v}{T}\tnormf{\nabla z}{\Nb^{\ell+1}(T)}\\
		&\leq C\Ct \Col  \beta^2 \alpha^{-1} \ell^{d/2} \exp(-\Cd \ell)\tnormf{\nabla v}{D}\tnormf{\nabla z}{D},
	\end{align*}
	with a constant $\Col>0$, which only depends on the regularity of the mesh $\calT_H$.
	The assertion now follows with the constant $\Cl = C\Ct \Col  \beta^2 \alpha^{-2}$, {which scales like $\Ct$ with respect to the polynomial degree $p$.}
\end{proof}

\end{document}